\documentclass[11pt,twoside]{article}
\usepackage{latexsym}
\usepackage{amssymb,amsbsy,amsmath,amsfonts,amssymb,amscd,amsthm}
\usepackage{epsfig, graphicx}
\usepackage{color}
\usepackage{enumerate}
\newcommand{\fer}[1]{(\ref{#1})}

\newcommand{\norm}[1]{\|#1 \|}
\newcommand{\R}{\mathbb{R}}

\newcommand{\Z}{\mathbb{Z}}
\newcommand{\Sph}{\mathbb{S}}

\newcommand {\e}  {\varepsilon}
\newcommand {\eps}  {\varepsilon}

\newcommand{\loc}{\mathrm{loc}}
\newcommand {\f}   {\frac}

\newcommand {\prt}{\partial_t }
\newcommand{\beq}{\begin{equation}}
\newcommand{\beqa}{\begin{eqnarray}}
\newcommand{\bea} {\begin{array}{ll}}
\newcommand{\beqan}{\begin{eqnarray*}}
\newcommand{\eeq}{\end{equation}}
\newcommand{\eeqa}{\end{eqnarray}}
\newcommand{\eeqan}{\end{eqnarray*}}
\newcommand{\eea} {\end{array}}
\newtheorem{theorem}{Theorem}[section]
\newtheorem{lemma}[theorem]{Lemma}

\newtheorem{proposition}[theorem]{Proposition}
\newtheorem{corollary}[theorem]{Corollary}
\newtheoremstyle{remarkb}
	{}
	{}
	{\normalfont}
	{}
	{\bfseries}
	{}
	{ }
	{}
\theoremstyle{remarkb}
\newtheorem{remark}[theorem]{Remark}
\numberwithin{equation}{section}
\newcommand{\cqfd}{{ \hfill
                       {\unskip\kern 6pt\penalty 500
                       \raise -2pt\hbox{\vrule\vbox to 6pt{\hrule width 6pt
                       \vfill\hrule}\vrule} \par}   }}
\title{\Large \bf Blow-up, concentration phenomenon and global existence  for the Keller-Segel model in high dimension}

\author{Vincent Calvez~$^{\text{a}}$, Lucilla Corrias~$^{\text{b}}$, Mohamed Abderrahman Ebde~$^{\text{c}}$
}
\date{\today}
\begin{document}
\maketitle
\pagestyle{plain}
\pagenumbering{arabic}

\begin{abstract}
This paper is devoted to the analysis of the classical Keller-Segel system over $\R^d$, $d\geq 3$. We describe as much as possible the dynamics of the system characterized by various criteria, both in the parabolic-elliptic case and in the fully parabolic case. The main results when dealing with the parabolic-elliptic case are: local existence without smallness assumption on the initial density, global existence under an improved smallness condition and comparison of blow-up criteria. A new concentration phenomenon criteria for the fully parabolic case is also given. The analysis is completed by a visualization tool based on the reduction of the parabolic-elliptic system to a finite-dimensional dynamical system of gradient flow type, sharing features similar to the infinite-dimensional system.
\end{abstract}



\section{Introduction}
\label{Sec:intro}

This paper aims to describe the dynamics of the Keller-Segel system in the whole space $\R^d$ and high dimension $d\ge3$ , with a particular emphasis given to blow-up and related facts. In its simplest and nondimensionalized formulation, the system reads as following
\begin{equation}\tag{KS}
\left\{\begin{array}{rcll}
\partial_t n& =& \Delta n - \nabla\cdot(n\nabla c)\, ,\medskip\\
\eps\, \partial_t c& =& \Delta c + n-\alpha c\, ,
\end{array}
\right.
\label{KS}
\end{equation}
and describes, at the macroscopic scale, a population of cells of density $n$ which attract themselves by secreting a diffusive chemical signal of concentration $c$. The nonnegative parameter $\eps$ is proportional to the ratio between the two diffusion coefficients of $n$ and $c$, appearing in the dimensionalized formulation of \eqref{KS}. It takes into account the different time scales of the two diffusion processes, $\eps = 0$ describing the chemical concentration evolution in a quasi-stationary approximation. In other words, $\eps$ is the parameter that influences the possible dynamics of $c$ and consequently of $n$. The chemical degradation rate $\alpha$ is also a nonnegative constant. It is related to the range of action of the signal, as we shall see later. 
 
In the quasi-stationary parabolic-elliptic case ($\eps = 0$), the chemical concentration should be understood as:
\begin{equation}
c=\left\{
\begin{array}{rcl}
E_d*n \,,\quad\alpha=0\,,  \medskip\\
B_d^\alpha*n\,,\quad\alpha>0\,,
\end{array}
\right.
\label{eq:chemical}
\end{equation}
where hereafter $*$ denotes the convolution with respect to the space variable and $E_d$ and $B_d^\alpha$ are respectively the Green's function for the Poisson's equation in $\R^d$, and the Bessel kernel:
\begin{align}
&E_d(x)=\mu_d\ \dfrac1{|x|^{d-2}}\,,\quad \quad \mu_d = \frac1{(d-2)|\Sph^{d-1}|}\,,
\label{eq:Ed} \medskip
 \\
&B_d^\alpha(x)= \int_0^{+\infty}\f1{(4\pi t)^{d/2}}\,  e^{-\frac{|x|^2}{4t} - \alpha t}\ dt\,.
\label{eq:Balpha}
\end{align}

The system is equipped with non negative initial data $n_0$ and $c_0$ (if $\eps>0$) and with fast decay conditions at infinity. Moreover, the initial cell density $n_0$ is supposed to be an integrable function, so that the total initial mass of cells $M$ is conserved along time:
\[ M =\int_{\R^d}n_0(x)\ dx=\int_{\R^d}n(x,t)\ dx\, ,\]
(see Theorem \ref{the:localexistence} for the complete set of assumptions).

The question of blow-up (or concentration) of solutions of the \eqref{KS} system has been a challenging issue in the field of mathematical biology since the seminal papers by Keller and Segel \cite{KS70,KS71}, as it describes the propensity of cells to aggregate when they interact through a chemical long-range signaling \cite{Nanjundiah, JL, NSY, Gajewski98}. The subtlety of this system has triggered a lot of theoretical works and many problems are only partially solved or still open. We refer to \cite{Horstmann03,P,HP09} for a complete overview of various mathematical results and modeling approaches.

It is well known, since the work by Childress and Percus \cite{CP81}, that the behaviour of the Keller-Segel system strongly depends on the space dimension. In dimension one, blow-up cannot occur for the solutions of \eqref{KS}, \cite{Nagai95}. In dimension two and for the parabolic-elliptic case ($\eps=0$) we have the so-called ``critical mass phenomenon'': when $\alpha = 0$ the $L^1-$norm of the cell density determines whether blow-up occurs or not, {\it i.e.} the cell density becomes unbounded in finite time if the mass is above some threshold (equals to $8\pi$), and it does not if the mass is below the same threshold. The same holds true when $\alpha>0$ up to an additional smallness condition on the second moment of $n_0$ for blow-up only, \cite{CC08}. This remarkable dichotomy is not yet proved to hold true when $\eps>0$. However, the mass still continue to be a key quantity. More precisely, when $\eps>0$ one can prove global existence if the mass is below $8\pi$ \cite{CC08} (using energy methods), but no clear result exists on the other side. In the case of a bounded domain, Horstmann proved that  blow-up occurs for some initial conditions having very negative energy  \cite{Horstmann02, Horstmann.Wang01}. To complicate the picture of the parabolic-parabolic case, the authors in \cite{BCD} showed the existence of positive forward self-similar solutions of \fer{KS} with $\alpha=0$, decaying to zero at infinity and having mass larger than $8\pi$ (which is no longer possible in the parabolic-elliptic case), see also \cite{NaitoSY} and the references therein.

In higher dimensions, the $L^1-$norm of $n_0$ is no longer a critical quantity. One available criterion for global existence involves the $L^{d/2}-$norm of the initial cell density \cite{{CP_CRAS06},CP08, CPZ04, KozoSugy}. On the other hand, blow-up is shown to occur when $\eps = 0$ if the second momentum of $n_0$ is small enough compared to the mass. Nothing is known about the possible blow-up of solutions when $\eps >0$, as in the two-dimensional case.

The main purpose of this paper is to shed a new light on the dynamics of the \eqref{KS} system  in both the parabolic-elliptic case and the parabolic-parabolic case, when $d\ge3$. We also propose a visualization tool to get the picture of the dynamics in the case $\eps = 0$. The tool is based on the reduction of \eqref{KS} firstly to a one-space  dimensional system and secondly to a finite-dimensional dynamical system sharing analogous features with \eqref{KS}.

More specifically, our main result concerning existence of solutions of \eqref{KS}  is the following local existence result for the parabolic-elliptic case, without smallness condition on the initial density $n_0$ and where the obtained weak solution is not to be intended in the integral sense (mild solution). For analogous results concerning mild solutions we refer to \cite{CP08,KozoSugy} and the references therein.

\begin{theorem}[Local existence]
Let $d\ge3$, $\alpha\ge0$ and $\eps=0$. Let $n_0$ be a nonnegative initial density in $(L^1\cap L^a)(\R^d)$, with $a>d/2$.   Assume in addition that $n_0\in L^1(\R^d,\psi(x) dx)$ where $\psi$ is a nonnegative function such that $\psi(x)\to+\infty$ uniformly as $|x|\to+\infty$, $e^{-\psi}\in L^1(\R^d)$ and $|\nabla \psi|\in L^\infty(\R^d)$. Then, there exists $T_{max}>0$ depending only on $\|n_0\|_{L^a(\R^d)}$ and a weak nonnegative solution in the distribution sense $(n,c)$ of \eqref{KS}, where $c$ is given by \eqref{eq:chemical}, such that
$$
n\in L^\infty((0,T_{max});(L^1\cap L^a)(\R^d))
\quad\hbox{and}\quad
n\in L^\infty_\loc((0,T_{max});L^p(\R^d))\,,\quad p\in(a,\infty)\,.
$$
Moreover, $n\in L^\infty((0,T_{max});L^1(\R^d,\psi(x) \,dx))$, according to the hypothesis on $n_0$. In addition, the total mass is conserved,
$$
M =\int_{\R^d}n_0(x)\ dx=\int_{\R^d}n(x,t)\ dx\,,\quad \hbox{a.e. } t\in(0,T_{max})\,,
$$
$nc$ and $n\log n$ belong to $L^\infty((0,T_{max});L^1(\R^d))$,   the energy $\mathcal E[n]\in L^\infty_\loc(0,T_{max})$, while the dissipation of energy $\int_{\R^d}n|\nabla(\log n-c)|^2\in L^1_\loc(0,T_{max})$. Finally, if  
$T_{max}<\infty$ then $\lim\limits_{t \nearrow T_{max}}\|n(t)\|_{L^a(\R^d)}=+\infty$.  
\label{the:localexistence}
\end{theorem}

The above local existence  becomes global if in addition we assume that the $L^{d/2}-$norm of the initial cell density is small enough. Our contribution in that direction is to obtain a smallness condition for $\|n_0\|_{L^{\f d2}(\R^d)}$ improved with respect to the existing one (see Corollary \ref{the:globalexistence} and Remark \ref{bestsmallnesscond}). Unfortunately, the results given in Theorem \ref{the:localexistence} and Corollary \ref{the:globalexistence} can not be extended to the fully parabolic case ($\eps>0$), essentially because in order to do so one has to control the temporal derivative of the chemical concentration $c$ (see Remark \ref{paraboliccase}). 

Concerning the blow-up issue, one of the key ingredient used here is the free energy naturally associated to \fer{KS}, in opposition with the two dimensional case where the free energy is fundamental for proving global existence results only  \cite{CC08,BDP}. The energy functional reads as follows
\begin{equation}
\mathcal E[n,c](t)=\int_{\R^d} n(x,t)\log n(x,t)\,dx - \int_{\R^d} n(x,t)c(x,t)\,dx +\, \f12\int_{\R^d} |\nabla c(x,t)|^2\,dx+ \f\alpha2 \int_{\R^d} c^2(x,t)\,dx\, ,
\label{eq:energy}
\end{equation}
and it satisfies the dissipation equation
\begin{equation}
\dfrac d{dt} \mathcal E[n,c](t) = - \int_{\R^d} n(x,t)\left|\nabla\left( \log n(x,t) - c(x,t) \right)\right|^2\, dx - \e \int_{\R^d} \left|\partial_t c(x,t)\right|^2\, dx\, .
\label{eq:energy dissipation} 
\end{equation}
The novelty in this paper, is that we will use also the corrected energy, defined as
\begin{equation}
\mathcal F[n,c](t)=
\log\left( \int_{\R^d} |x|^2 n(x,t)\, dx\right) +\f{2}{dM}\mathcal{E}[n,c]\,.
\label{eq:corrected energy} 
\end{equation}
 Under the quasi-stationary assumption $\eps = 0$, the free energy \eqref{eq:energy} reduces to the difference between  the entropy and the potential energy:
\begin{equation}
\mathcal E[n](t)=\int_{\R^d} n(x,t)\log n(x,t)\,dx - \dfrac12\int_{\R^d} n(x,t)c(x,t)\,dx\, .
\label{eq:energyPE}
\end{equation}

The blow-up results analyzed here with the help of the above tools can be summarized as follows. In the parabolic-elliptic case ($\eps=0$ and $\alpha\ge0$), blow-up occurs in finite time if one of the following criteria is fullfilled:
$$
\int_{\R^d} |x|^2 n_0(x)\, dx < K_1^\alpha(d,M) M^{\frac{d}{d-2}}\,, 
$$
or
$$
\int_{\R^d} |x|^2 n_0(x)\, dx < K_2(d)\, M^{\,1+\frac{2}d}\, \exp\left(-\dfrac2{dM}\,\mathcal{E}[n_0]\right)\,,
$$
where the constants $K_1^\alpha(d,M)$ and $K_2(d)$ are defined in Propositions \ref{prop:1BUcriterion} and \ref{prop:2BUcriterion} respectively. After having derived both criteria, we show that none of these two criteria contain the other. However we have some evidence that the second one (involving the free energy) appears to be better and the natural one in some sense. Indeed, this is clearly the case when dealing with the discrete model in Section \ref{Sec:num}. Moreover, it is possible to extend  the second criteria to the parabolic-parabolic case ($\eps>0$ and $\alpha=0$ for the sake of simplicity). Doing that, we are able to prove the following concentration result for the critical norm of the cell density (recovering the parabolic-elliptic case as $\eps\searrow 0$).

\begin{theorem}[Parabolic concentration]
Let $\e>0$, $\alpha=0$ and $d\ge3$. Assume that the nonnegative initial densities $(n_0,c_0)$ have finite energy $\mathcal E[n_0,c_0]$ and satisfy
\begin{equation}
\int_{\R^d} |x|^2 n_0(x)\, dx< K_2(d)\
M^{1+ \f{2}d}\exp{\left(-\f2{dM}\mathcal E[n_0,c_0]\right)}
\exp\left(-\varepsilon^\gamma\right)\, ,  
\label{eq:concentration_condition}  
\end{equation}
where the constant $K_2(d)$ is given in \eqref{2BU} and $\gamma\in(0,1)$. Let $(n,c)$ be a sufficently smooth solution of \eqref{KS} generated by $(n_0,c_0)$ and $T_{max}$ the maximal time of existence (possibly infinite). Then there exists a constant $C(d)>0$ such that 
\begin{equation}
\sup_{t\in[0,T_{max})}\,\|n(t)\|_{L^{\f d2}(\R^d)}
\geq \f{2(d-2)}{C^2(d)}\,\e^{\gamma-1}\,.
\label{eq:concentration}  
\end{equation}
\label{th:concentration}
\end{theorem}

The plan of the paper is as follows. We first list in Section~\ref{Sec:inequalities} some useful  sharp functional inequalities, which appear to be key tools for  the quantification of criteria involving critical quantities. In Section~\ref{Sec:LGexistence} we present the new local in time existence result for the parabolic-elliptic Keller-Segel system without any smallness condition on the initial data (Theorem \ref{the:localexistence}) and we review conditions for global existence (Corollary \ref{the:globalexistence}). In Section \ref{Sec:Blow-up} we derive the above two blow-up criteria for the parabolic-elliptic case (Propositions \ref{prop:1BUcriterion} and \ref{prop:2BUcriterion}) and we discuss the complementarity of those criteria (Section \ref{sec:complementarity}). In Section \ref{Sec:Blow-up para-para} we extend the criterion involving the free energy to the parabolic-parabolic case into a concentration result (Theorem \ref{th:concentration}). Finally, in Section \ref{Sec:num}, we shall analyze  finite-dimensional energy functionals whose gradient flow possesses a structure similar to the parabolic-elliptic \eqref{KS} system. We derive accordingly one criterion ensuring global existence and two criteria for blow-up. The overall dynamical picture is drawn in Fig. \ref{fig:BU non homogene}.

\section{Functional inequalities and preliminaries}
\label{Sec:inequalities}
When analyzing the Keller-Segel system \fer{KS} over $\R^2$ (whatever is $\eps\ge0$), the free energy \fer{eq:energy} is the key tool to prove the global existence of \emph{weak  free energy solutions} (see \cite{BCM} for the definition). Indeed, dual functional inequalities (the logarithmic Hardy-Littlewood-Sobolev inequality and the Moser-Trudinger-Onofri inequality) applied to the energy, together with the time decreasing behaviour of the energy itself, give the key a priori estimates on the entropy term $\int_{\R^2}n\log n$, the potential interaction term $\int_{\R^2}nc$ and on the $H^1$ norm of $c$, \cite{CC08}. 

Conversely, when the \fer{KS} system is analyzed over $\R^d$, $d\ge3$, these inequalities are no longer adapted to the energy functional. Therefore, one has to establish a priori estimates on some $L^p$ norms of the solutions in order to get existence results of more classical \emph{weak solutions in the distribution sense}. 

We shall show in this section that, even in the high dimensional case, when $\eps=0$ sharp functional inequalities, corresponding in some sense to those mentioned above for the two dimensional case, are the key tools for both the local and the global existence results. More precisely, our purpose is to use optimally the sharpness of these inequalities and their ``best constant'' in order to quantify possible thresholds of eventually critical norms of the cell density $n$.  Consequently, the results given here for the parabolic-elliptic Keller-Segel system  ($\eps = 0$) are not only refinement of known results but also new ones. This is the case for instance for the local existence of weak solutions without any smallness condition on the initial density $n_0$ and the hypercontractivity property of the solutions. On the other hand, the already known smallness condition on $n_0$ sufficient for the global existence will be improved thanks to a more careful estimate and conjectured to be the optimal one (see Remark \ref{bestsmallnesscond}).

Concerning the parabolic-parabolic \fer{KS} system ($\eps>0$), it appears that the estimates valid for the parabolic-elliptic case can not be reproduced optimally. Therefore, we are not able to obtain here a local existence result without any smallness conditions, nor to improve the known smallness condition for the global existence when $\eps>0$ (see Remark \ref{paraboliccase}). This is one major difference with the two dimensional case, where the energy \fer{eq:energy} can be handle in the same way, whatever is $\eps\ge0$ \cite{CC08}.

The mentioned inequalities used here and throughout this paper are firstly the classical Sobolev inequality 
\begin{equation}
\norm{f}_{L^{\frac{2d}{d-2}}(\R^d)}\le C_S(d)\norm{\nabla f}_{L^{2}(\R^d)}\,,\quad\quad
C^2_S(d):=\f4{d(d-2)|\Sph^d|^{2/d}}\,,
\label{in:Sob}
\end{equation}
and the following special case of the sharp Hardy-Littlewood-Sobolev inequality \cite{LiebLoss}
\begin{equation}
\left|\iint_{\R^d\times\R^d}f(x)|x-y|^{-\lambda}g(y)\ dxdy\right|\le C_{HLS}(d,\lambda)\|f\|_{L^p(\R^d)}\|g\|_{L^p(\R^d)}\, ,\quad p=\frac{2d}{2d-\lambda}\,,\quad 0<\lambda<d\, ,
\label{in:HLS}
\end{equation}
where the best constant $C_{HLS}(d,\lambda)$ has been obtained by Lieb \cite{Lieb}. 
%
More specifically, because of \fer{eq:Ed} and \eqref{eq:Balpha}, inequality \fer{in:HLS} will be used in the sequel with $\lambda=d-2$. In this case \fer{in:Sob} and \fer{in:HLS} turn to be dual, as proved for instance in \cite{LiebLoss} Theorem 8.3, with
\begin{equation}
C_{HLS}(d,d-2)=\pi^{\f d2-1}\Gamma^{-1}\left(\f d2+1\right)\left(\f{\Gamma(d)}{\Gamma(d/2)}\right)^{\f 2d}\ ,
\label{eq:CHLSd-2}
\end{equation}
so that the following relation can be established  
\begin{equation}
C_S^2(d)=\mu_d\,C_{HLS}(d,d-2)\ .
\label{eq:CS-HLS}
\end{equation}
 Let us observe that when $\lambda\searrow0$, \fer{in:HLS} boils down to the logarithmic Hardy-Littlewood-Sobolev inequality \cite{CarlenLoss,Beckner} used in the two dimensional case \cite{BDP,CC08}.

The other functional inequality we shall use is simply the following Gagliardo-Nirenberg interpolation inequality, 
\begin{equation}
\norm{v}_{L^{\f{2(p+1)}p}(\R^d)}\le C_{GN}(p,d)\|\nabla v\|_{L^2(\R^d)}^{\f p{p+1}}\|v\|_{L^{\f dp}(\R^d)}^{\f1{p+1}}\, ,\quad 1\le p\le d\, .
\label{GN}
\end{equation}
We will focus on the particular case $p = \f d2$ in \eqref{GN} and pay much attention to the  best constant $C_{GN}(d/2,d)$. Generally, it is not possible to compute the optimal constants in the family of inequalities \fer{GN}. However and very fortunately, it is possible in the special case $p = \f d2$ under interest here, since  it corresponds to the  Gagliardo-Nirenberg interpolation inequality 
\begin{equation}
\norm{v}_{L^{2\sigma+2}(\R^d)}^{2\sigma+2}\le C_{\sigma,\,d}^{2\sigma+2}
\|\nabla v\|_{L^2(\R^d)}^{\sigma d}\|v\|_{L^2(\R^d)}^{2+\sigma(2-d)}\ ,
\quad 0<\sigma<\f2{d-2}\,,\quad d\ge2\,,
\label{GN2}
\end{equation}
for $\sigma=\f 2d$. Weinstein has characterized the sharp constant in this context \cite{W}, so that
\begin{equation}
C_{GN}\left(\f d2,d\right)=\ C_{\f2d\,,\,d}\ =\left(1+\f2d\right)^{\f d{2(d+2)}}\|\psi\|_{L^2(\R^d)}^{-\f2{d+2}}\,,
\label{GNcost}
\end{equation}
where $\psi$ is a $H^1(\R^d)\cap C^\infty(\R^d)$ positive and radial function, solution of minimal $L^2$ norm of the equation $\Delta\psi-\f2d\psi+\psi^{\f4d+1}=0$.
\\

Next come two interpolation lemmas, useful in the sequel for the control of both the entropy $\int_{\R^d}n\log n$ and the potential $\int_{\R^d}nc$ in the free energy. 
\begin{lemma}[Entropy lower bound]
Let $f$ be any nonnegative $L^1(\R^d)$ function such that  \break $I~=~\int_{\R^d}|x|^2f(x)\, dx<\infty$ and $\int_{\R^d}f \log f <\infty$. Let $M=\int_{\R^d}f $. Then,
\begin{equation}
\int_{\R^d}f \log f  +\delta I\ge M\log M+\f{dM}2\log(\delta\,\pi^{-1})\,,
\label{est:entropy1}
\end{equation}
for any $\delta>0$, and
\begin{equation}
\int_{\R^d}f \log f +\f{dM}2[\log I+1]\ge M\log M+\f{dM}2\log\left(\f{dM}{2\pi}\right)\ .
\label{est:entropy2}
\end{equation}
\label{est:entropy}
\end{lemma}
\begin{proof}
Applying the Jensen's inequality with the probability density $\mu(x)=\delta^{d/2}\pi^{-  d/2}e^{-\delta|x|^2}$, one obtain:
\begin{align*}
\int_{\R^d}f \log f  &=\int_{\R^d}\f {f }{\mu }\log \left(\f {f }{\mu }\right)\mu +\int_{\R^d}f \log\mu  \\ & \ge
M\log M+\int_{\R^d}f \log\mu \, ,
\end{align*}
and \fer{est:entropy1} follows. The optimization of \fer{est:entropy1}  with respect to $\delta >0 $ under the fixed constrains $I$ and $M$, yields \fer{est:entropy2}, the optimal $\delta$ being $\delta=\f{dM}{2I}$.
\end{proof}

\begin{remark}
It is worth noticing that the entropy lower bound \eqref{est:entropy1}, and so \eqref{est:entropy2} too, is actually a sharp inequality since for any fixed $\delta>0$ equality holds true in \eqref{est:entropy1} if (and only if) $f=M\mu$. 
\label{rem:sharpentropylowerbound}
\end{remark}

\begin{lemma}[Potential confinement]
Let $f$ be any nonnegative  function such that $f\in (L^1\cap L^{\f d2})(\R^d)$ and  $I=\int_{\R^d}|x|^2f(x)\, dx<\infty$. Let $M=\int_{\R^d}f(x)\, dx$. Then,
\begin{equation}
2^{\,1-\f d2}\,M^{\,\f d2+1}\,I^{\,1-\f d2}\le\,\iint_{\R^d\times\R^d}f(x)\f1{|x-y|^{d-2}}f(y)\,dxdy\le C_{HLS}(d,d-2)\,M\,\|f\|_{L^{\f d2}(\R^d)}\ .
\label{est:potential}
\end{equation}
\label{lem:estpotential}
\end{lemma}
\begin{proof}
The right inequality is a direct consequence of the Hardy-Littlewood-Sobolev inequality \fer{in:HLS} and of standard interpolation (see also \cite{BCL}). For the left inequality we have simply \cite{B95}
\begin{align*}
M^2 & = \iint_{\R^d\times\R^d}f(x)f(y)\, dxdy \\
& \leq \left(\iint_{\R^d\times\R^d} f(x) |x-y|^2 f(y)\, dxdy \right)^{1-2/d} \left(\iint_{\R^d\times\R^d} f(x) \dfrac1{|x-y|^{d-2}} f(y)\, dxdy  \right)^{2/d}\\
&=\left(2IM -2\,\left|\int_{\R^d} xf(x)\,dx\right|^2
\right)^{1-2/d} \left( \iint_{\R^d\times\R^d} f(x) \dfrac1{|x-y|^{d-2}} f(y)\, dxdy \right)^{2/d}
\end{align*}
and the proof follows.
\end{proof}

We conclude this section by recalling some  few basic properties satisfied by the Bessel kernel $B_d^\alpha$ defined in \eqref{eq:Balpha} and usefull in the sequel. Probably they are classical ones. However, since we do not found any references for them, we list and prove the required properties in the lemma below.
\bigskip
\begin{lemma}[Properties of the Bessel kernel]
The following relations for $B_d^\alpha$ and $\nabla B_d^\alpha$ hold true for $\alpha\ge0$ in any space dimension $d\ge3$.
\begin{enumerate}[(i)]
\item Expansion of $B_d^\alpha$ with respect to $E_d$: 
\begin{equation}
B_d^\alpha = E_d - \alpha\, E_d*B_d^\alpha\quad\quad a.e.\,,
\label{eq:Expansion alpha}
\end{equation}
\item Gradient formula:
\begin{equation}
\nabla B_d^\alpha(x)=-\f1{|\Sph^{d-1}|}\,\f x{|x|^d}\,g_\alpha(|x|)\,,
\label{eq:gradBalpha}
\end{equation}
where $g_\alpha(|x|):=\Gamma (d/2)^{-1}\int_0^{+\infty}s^{\f d2-1}e^{-s-\alpha\f{|x|^2}{4s}}\,ds$\,,  is a positive radial function exponentially decreasing from 1 to 0 as $|x|\to\infty$.
\item Corrected Euler's homogeneous function theorem:
\begin{equation} x\cdot \nabla B_d^\alpha(x) = -(d-2) B_d^\alpha(x) - 2 \alpha\, (B_d^\alpha*B_d^\alpha)(x)\, . 
\label{eq:Euler alpha}
\end{equation}
\end{enumerate}
\label{lem:Besselproperties}
\end{lemma}

\begin{proof} First of all, let us observe that the convolutions in \eqref{eq:Expansion alpha} and \eqref{eq:Euler alpha} are well defined since both $B_d^\alpha$ and $E_d$ belong to $L^{\f d{d-2}}_w(\R^d)$, (the weak $L^{\f d{d-2}}$ space), and $B_d^\alpha\in L^p(\R^d)$ for $1\le p<\f d{d-2}$, (see \cite{LiebLoss}). Then, identity \eqref{eq:Expansion alpha} make sense a.e.\,. Next, identity \eqref{eq:Expansion alpha} is an immediate consequence of the facts that (i) $E_d$ and $B^\alpha_d$ are the Green's functions for $-\Delta$ and $(-\Delta+\alpha)$, respectively, and (ii) $B_d^\alpha*v$ is the unique solution of $(-\Delta+\alpha)u=v$ in ${\cal D}'(\R^d)$, belonging to $L^r(\R^d)$ for some $r\ge1$, (see \cite{LiebLoss}).

Formula \eqref{eq:gradBalpha} is a straightforward computation. To prove \eqref{eq:Euler alpha}, first notice that the Fourier transform of $B_d^\alpha$ is $\left(\alpha + 4\pi^2|\xi|^2\right)^{-1}$, when defining the Fourier transform as $\hat f(\xi)=\int_{\R^d}e^{-2\pi i\,x\cdot\xi}f(x)\,dx$,  (see \cite{LiebLoss}). Then, denoting $f^\vee$ the inverse Fourier transform, we have in the Fourier variable $\xi$,
\begin{align*}
(x\cdot\nabla_x B_d^\alpha)\, \widehat\ \,(\xi) & = 
(x\cdot\nabla_x((\alpha + 4\pi^2|\xi|^2)^{-1})^{\vee})\,\widehat\ \,(\xi)
= 2\pi\,i\,(x\cdot(\xi(\alpha + 4\pi^2|\xi|^2)^{-1})^{\vee})\,\widehat\ \,(\xi)
\\
&=- \nabla_\xi\cdot\left(\f\xi{\alpha + 4\pi^2|\xi|^2}\right)
= - \dfrac d{\alpha + 4\pi^2|\xi|^2} + \f{8\pi^2|\xi|^2}{\left(\alpha + 4\pi^2|\xi|^2\right)^2} \\ 
& = -\dfrac{d-2}{\alpha + 4\pi^2|\xi|^2} - \f {2\alpha} {\left(\alpha + 4\pi^2|\xi|^2\right)^2}
\end{align*}
and identity \eqref{eq:Euler alpha} is proved.

\end{proof}

\begin{remark} Actually, in order to prove Lemma \ref{lem:Besselproperties}, it is possible to develop an alternative argument based on the homogeneity property of $B_d^\alpha$, namely: $B_d^\alpha(\lambda x)=\lambda^{2-d}B_d^{\alpha \lambda^2}(x)$. From this point of view, formula \eqref{eq:Euler alpha} can be seen as a corrected Euler's homogeneous function theorem. Obviously \eqref{eq:Euler alpha} reduces exactly to the Euler's homogeneous function theorem for $E_d$ as $\alpha = 0$, $E_d$ being homogeneous of degree $(2-d)$. 
\end{remark}

\section{Local and global existence for the parabolic-elliptic system ($\eps=0$)}
\label{Sec:LGexistence}
This section is entirely devoted to the proof of Theorem \ref{the:localexistence}. The proof is given in several steps, beginning with a cascade of a priori estimates. The functional inequalities given in the previous section and the mass conservation will play a fundamental role. The rate of decay at infinity of the initial density $n_0$ is also essential. Let us observe that one can weaken the usual hypothesis $n_0\in L^1(\R^d,|x|^2dx)$, as we have done here, and obtain existence result also with initial density having infinite second moment but logarithmic decay (see also Subsection \ref{subsection:entopy estimates} and Remark~\ref{rm:logdecay}). This is absolutely not the case when dealing with blow-up results, always obtained under the hypothesis of finite initial second moment.

As a byproduct of the a priori estimates, we obtain the following global existence result under an improved smallness condition on the initial density. Since the corollary is a straightforward consequence of Theorem \ref{the:localexistence}, the proof will be  omitted.
\begin{corollary}[Global existence]
Under the same hypothesis on the initial density $n_0$ as in Theorem~\ref{the:localexistence}, if in addition
\begin{equation}
\norm{n_0}_{L^{\frac d 2}(\R^d)}<\,\f{8}{d}\,C_{GN}^{-2(1+\f2d)}\left(\f d2,\,d\right)\ ,
\label{n_0condGN}
\end{equation}
the above weak solution exists globally in time. Moreover, the $L^{\f d2}$-norm of $n$ is  time decreasing.
\label{the:globalexistence}
\end{corollary}
\begin{remark}
In high dimension ($d\ge6$), Corollary \ref{the:globalexistence} can be proved under the sharp hypothesis on the initial density $n_0\in(L^1\cap L^{\f d2})(\R^d)$, as it follows by an attentive analysis of the proof of Theorem~\ref{the:localexistence}, see subsection \ref{prooflocalexistence}.
\end{remark}
%

\subsection{Local in time estimate of $\norm{n}_{L^p(\R^d)}$}
\label{subsec:localexistence}
Since $-\Delta c=n-\alpha c$ , we have for any nonnegative solutions $n$
\begin{equation}
\f{d}{dt}\int_{\R^d} n^p=- 4\f{p-1}{p}\int_{\R^d}|\nabla n^{\f p2}|^2 +(p-1)\int_{\R^d} n^{p+1}-\alpha(p-1)\int_{\R^d}n^p\,c \ .
\label{est:np}
\end{equation}
Next, using standard interpolation and the Sobolev inequality \fer{in:Sob}, we get 
\begin{equation}
\int_{\R^d} n^{p+1}=\int_{\R^d}(n^{\f p2})^{\f2p(p+1)}
\le\norm{n^{\f p2}}^{(p+1-\f d2)\f2p}_{L^2(\R^d)}\
\norm{n^{\f p2}}_{L^{\f{2d}{d-2}}(\R^d)}^\f dp
\le\norm{n}_{L^p(\R^d)}^{(p+1-\f d2)}
\left(C_S(d)\norm{\nabla n^{\f p2}}_{L^{2}(\R^d)}\right)^{\f dp}\ ,
\label{est:n(p+1)}
\end{equation}
for any $\max\{\f d2-1\,;\,1\}\le p<\infty$. Inequality \fer{est:n(p+1)} gives us for $p>\f d2$
\begin{equation}
\int_{\R^d} n^{p+1}\le
\f{\delta^{r'}}{r'}\norm{n}_{L^p(\R^d)}^{(p+1-\f d2)r'}+
\f1{r\delta^r}\, C^2_S(d)\norm{\nabla n^{\f p2}}_{L^{2}(\R^d)}^2\ ,
\label{est:d2+}
\end{equation}
where $r=\f{2p}d$, $r'=\f r{r-1}$ and $\delta>0$. Plugging estimate \fer{est:d2+} into \fer{est:np}, one obtains 
\begin{equation}
\f{d}{dt}\int_{\R^d} n^p\le\f{(p-1)}p\Big[\f{dC^2_S(d)}{2\delta^r} -4\Big]\norm{\nabla n^{\f p2}}_{L^{2}(\R^d)}^2+
(p-1)\f{\delta^{r'}}{r'}\left(\int_{\R^d} n^p\right)^{1+(p-\f d2)^{-1}}\ .
\label{est:nd/2+}
\end{equation}
Then, for any fixed $p>\f d2$,  it is enough to choose $\delta=\delta(p)>0$ such that $\f{dC^2_S(d)}{2\delta^r} -4\le0$ in order to have from \fer{est:nd/2+} a local in time boundeness of the $L^p$ norm of $n$, $p>\f d2$, whenever $n_0\in L^p(\R^d)$, i.e.
\begin{equation}
\int_{\R^d} n^p(t)\le\,h_p(t)\int_{\R^d} n_0^p
\ ,
\label{est:nd/2+bis}
\end{equation}
for $t\in[0,T_p)$, where 
\begin{equation}
h_p(t):=\left[1-\left(1-\f1p\right)\delta^{r'}\,t\left(\int_{\R^d} n_0^p\right)^{(p-\f d2)^{-1}}\right]^{(\f d2-p)}
\quad\hbox{and}\quad
T_p:=p(p-1)^{-1}\delta^{-r'}\left(\int_{\R^d} n_0^p\right)^{(\f d2-p)^{-1}}\ .
\label{eq:Ta}
\end{equation}
This is obviously only a local in time estimate since $h_p(t)\to+\infty$ as $t\nearrow T_p$ .

\subsection{Global in time estimate of $\norm{n}_{L^p(\R^d)}$ under smallness condition}
\label{Global}
In order to obtain a global in time estimate of $\norm{n}_{L^p(\R^d)}$, one has to be more careful in the estimate of the $L^{p+1}$ norm of $n$ than in \fer{est:n(p+1)} and use directly the sharp Gagliardo-Nirenberg inequality \fer{GN}. Indeed, we have
\begin{equation}
\int_{\R^d} n^{p+1}=\int_{\R^d}(n^{\f p2})^{\f2p(p+1)}
\le C_{GN}^{\f2p(p+1)}(p,d)\,\norm{\nabla n^{\f p2}}_{L^{2}(\R^d)}^2\,\norm{n}_{L^{\frac d 2}(\R^d)}\ .
\label{est:n(p+1)bis}
\end{equation}
Plugging \fer{est:n(p+1)bis} into \fer{est:np} we get for any $1\le p\le d$
\begin{equation}
\f{d}{dt}\int_{\R^d} n^p(t)\le(p-1)\Big[C_{GN}^{\f2p(p+1)}(p,d) \norm{n(t)}_{L^{\frac d 2}(\R^d)}-\f4p\Big]\norm{\nabla n^{\f p2}}_{L^{2}(\R^d)}^2\ .
\label{est:nd/2}
\end{equation}
As a consequence, choosing $p=\f d2$ into \fer{est:nd/2}, whenever the initial density satisfies \fer{n_0condGN}, 
the $\norm{n(t)}_{L^{\frac d 2}(\R^d)}$ norm decreases for all times $t\geq 0$.

\subsection{Local in time hyper-contractivity property}
Let $n_0\in (L^1\cap L^a)(\R^d)$, with $a>\f d2$ arbitrarily closed to $\f d2$. From subsection \ref{subsec:localexistence}, there exists a finite time $T_a>0$, defined in \eqref{eq:Ta}, and a local solution $n\in L^\infty((0,T_{a});(L^1\cap L^a)(\R^d))$. We are going to prove that, for any $p\in(a,+\infty)$, there exists a constant $C$ not depending on $\norm{n_0}_{L^p(\R^d)}$ such that
\begin{equation}
\int_{\R^d}n^p(t)\le C(1+t^{\,1-p}\,)\ ,\quad\quad\hbox{a.e. }\ t\in(0,T_a)\,.
\label{est:hyper}
\end{equation}

 Let $\delta\in(0,T_a)$ be arbitrarily small. Owing to the local in time boundedness of the $L^a$-norm of $n$, there exists a modulus of  ``\,$\f d2$-equintegrability'' $\omega(K;T_a-\delta)$ such that for $K\ge1$ it holds 
\begin{equation}
\sup\limits_{0\le t\le (T_a-\delta)}\,\norm{(n-K)_+(t)}_{L^{\f d2}(\R^d)}^{\f d2}\le \omega(K;T_a-\delta)
\quad\hbox{and}\quad
\lim\limits_{K\to+\infty}\omega(K;T_a-\delta)=0\ .
\label{equi}
\end{equation}
Indeed, for $t\in[0,T_a-\delta]$, using the increasing behaviour of $h_a(t)$ defined in \fer{eq:Ta}, we have
$$
\norm{(n-K)_+(t)}_{L^{\f d2}(\R^d)}^{\f d2}\le\f1{K^{a-\f d2}}\int_{\R^d}n^a(t)\le
\f{h_a(T_a-\delta)}{K^{a-\f d2}}\int_{\R^d}n_0^a\,=:\omega(K;T_a-\delta)\ .
$$

Next, let us estimate the time evolution of $\norm{(n-K)_+}_{L^p(\R^d)}$, from the differential inequality
\begin{equation}
\begin{array}{rcl}
\f{d}{dt}\int_{\R^d}(n-K)_+^p\le
-4\f{(p-1)}p\int_{\R^d}|\nabla(n-K)_+^{\f p2}|^2
+(p-1)\int_{\R^d}(n-K)_+^{p+1}\\ \\
+K(2p-1)\int_{\R^d}(n-K)_+^{p}+
pK^2\int_{\R^d}(n-K)_+^{p-1}\ .
\end{array}
\label{est:(n-K)p}
\end{equation}
For any fixed $p>\max\{2;a\}$, interpolating the $L^p$ and $L^{p-1}$ norms of $(n-K)_+$ between $L^1$ and $L^{p+1}$, one get easily from \fer{est:(n-K)p}
\begin{equation}
\f{d}{dt}\int_{\R^d}(n-K)_+^p\le
-4\f{(p-1)}p\int_{\R^d}|\nabla(n-K)_+^{\f p2}|^2
+4(p-1)\int_{\R^d}(n-K)_+^{p+1}
+4K^pM\ .
\label{est:(n-K)pbis}
\end{equation}
It remains to take advantage of the negative term in the right hand side of \fer{est:(n-K)pbis}. From standard interpolation and the Sobolev inequality \fer{in:Sob} we have for any $p\ge\max\{\f d2-1\,;\,1\}$ 
\begin{equation}
\int_{\R^d} v^{p+1}\le\norm{v}^p_{L^{\frac{pd} {d-2}}(\R^d)}\
\norm{v}_{L^{\frac d 2}(\R^d)}=
\norm{v^{\f p2}}^2_{L^{\frac{2d}{d-2}}(\R^d)}\
\norm{v}_{L^{\frac d 2}(\R^d)}
\le C^2_S(d)\norm{\nabla v^{\f p2}}_{L^{2}(\R^d)}^2\ \norm{v}_{L^{\frac d 2}(\R^d)}\ .
\label{est:n(p+1)Sob}
\end{equation}
Let us observe that here it is more convinient to use \eqref{est:n(p+1)Sob} than \eqref{est:n(p+1)bis}, in order to have a larger range of $p$ index for which the inequality holds true. 

Using \fer{est:n(p+1)Sob} and \fer{equi}, we obtain for $t\in[0,T_a-\delta]$ the estimate
$$
\int_{\R^d}(n-K)_+^{p+1}\le C_S^2(d)\norm{\nabla(n-K)_+^{\f p2}}_{L^2(\R^d)}^2\,\omega^{\f2d}(K;T_a-\delta)\ ,
$$
so that \fer{est:(n-K)pbis} becomes
\begin{equation}
\f{d}{dt}\int_{\R^d}(n-K)_+^p\le
-4(p-1)\left[p^{-1}C_S^{-2}(d)\omega^{-\f2d}(K;T_a-\delta)-1\right]
\int_{\R^d}(n-K)_+^{p+1}
+4K^pM\ .
\label{est:(n-K)pter}
\end{equation}
Taking $K$ sufficiently large, the quantity $\eta:=\left[p^{-1}C_S^{-2}(d)\omega^{-\f2d}(K;T_a-\delta)-1\right]$ is positive and using once again the interpolation of the $L^p$ norm of $(n-K)_+$ between $L^1$ and $L^{p+1}$, we obtain
\begin{equation}
\f{d}{dt}\int_{\R^d}(n-K)_+^p\le
-4(p-1)\,\eta\,M^{-\f1{p-1}}\left(\int_{\R^d}(n-K)_+^p\right)^{\f p{p-1}}
+4K^pM\ .
\label{est:(n-K)pquater}
\end{equation}
We are finally able to prove that there exists a positive finite constant $C=C(T_a,K,M,p)$, not depending on $\int_{\R^d}(n_0-K)_+^p$, such that for any $p\ge\max\{2;a\}$ and  for $t\in(0,T_a-\delta]$
$$
\int_{\R^d}(n-K)_+^p(t)\le \,\f C{t^{\ p-1}}\ ,
$$
simply by comparison of positive solutions of \fer{est:(n-K)pquater} with positive solutions of the differential equation $u'(t)+4(p-1)\,\eta\,M^{-\f1{p-1}}\,u^{\f p{p-1}}(t)=4K^pM$. Consequently and as usual, the hypercontractivity estimate \fer{est:hyper} holds true for any $p>\max\{2;a\}$.  For $p\in(a,\max\{2;a\}]$ it follows by interpolation.

%
\subsection{A priori estimates for the chemical density}
Let $n_0\in (L^1\cap L^a)(\R^d)$. Then, from the previous estimates, any corresponding nonnegative solutions $n$ belongs to any $L^p(\R^d)$, $p\in(1,\infty)$, for a.e. $t\in(0,T_a)$. Consequently, taking $p\in(1,\f d2)$, $c$ in \fer{eq:chemical} is well defined for any $\alpha\ge0$, it belongs to $L^q(\R^d)$, with $q=\f{pd}{d-2p}$\,, and
$$
\|c(t)\|_{L^q(\R^d)}\le C(d,\alpha)\,\|n(t)\|_{L^p(\R^d)}\,,\quad\hbox{a.e. } t\in(0,T_a)\ ,
$$
by the weak Young inequality \cite{LiebLoss}, (see also Lemma \ref{lem:Besselproperties}). In the same way, $\nabla c$ is well defined in the distributional sense as
\begin{equation}
\nabla c=\left\{
\begin{array}{rcl}
\nabla E_d*n\,,\quad\alpha=0\,,  \medskip\\
\nabla B_d^\alpha*n\,,\quad\alpha>0\,.
\end{array}
\right.
\nonumber
\end{equation}
and $|\nabla c|\in L^r(\R^d)$, with $r=\f{pd}{d-p}$, for $p\in(1,d)$: 
\begin{equation}
\|\nabla c(t)\|_{L^r(\R^d)}\le C(d,\alpha)\,\|n(t)\|_{L^p(\R^d)}\,,\quad\hbox{a.e. } t\in(0,T_a)\, .
\label{est:gradc}
\end{equation}
%

\subsection{Entropy, potential, energy and energy dissipation estimates}
\label{subsection:entopy estimates}
To conclude the a priori estimates, we will show how a local in time control on the $L^{\f d2}$ norm of $n$ together with a control on the decay of $n_0$ as $|x|\to\infty$, give local in time estimates on the entropy, the potential, the energy and the dissipation of the  energy. 

  Recall that under the assumptions of Theorem \ref{the:localexistence} we have that the initial density satisfies $n_0\in  (L^1\cap L^{\f d2})(\R^d)\cap L^1(\R^d,\psi(x)dx)$, where $\psi$ is nonnegative and satisfies $e^{-\psi}\in L^1(\R^d)$ and $|\nabla \psi|\in L^\infty(\R^d)$. Let $n$ be a nonnegative solution such that $n(t)\in (L^1\cap L^{\f d2})(\R^d)$ for a.e. $t\in(0,T)$, $T>0$. Then, the local in time bound on the potential follows directly from the definition of $E_d$ and \eqref{eq:Expansion alpha} (for the case $\alpha>0$) and from the potential confinement Lemma \ref{lem:estpotential} that gives
\begin{equation}
0\le\int_{\R^d}n(t)\,c(tt) \le\mu_d\,C_{HLS}(d,d-2)M\|n(t)\|_{L^{\f d2}(\R^d)}\,.
\label{est:nc}
\end{equation}

Next, the positive and negative contributions of the inital entropy $\int_{\R^d}n_0\log n_0$ are bounded. Indeed, for the positive contribution it simply holds true, for any $d\ge3$,   that
\begin{equation}
0\le\int_{\R^d}(n_0\log n_0)_+\le\|n_0\|_{L^\f d2(\R^d)}\,.
\label{est:nlognpos}
\end{equation}
Moreover, setting $v=n_0\,1\!\!1_{\{n_0\le1\}}$ and $m=\int_{\R^d}v$, the Jensen inequality gives us
$$
\int_{\R^d}(v\log v-v\psi)=\int_{\R^d}v\log\left(\f v{e^{-\psi}}\right)\ge m\log m-m\log(\|e^{-\psi}\|_{L^1(\R^d)}),
$$
which implies
\begin{equation}
0\le\int_{\R^d}(n_0\log n_0)_-\le m\log(\|e^{-\psi}\|_{L^1(\R^d)})-m\log m
-\int_{\R^d}n_0\,1\!\!1_{\{n_0\le1\}}\psi
\le C(m,\|e^{-\psi}\|_{L^1(\R^d)})+\int_{\R^d}n_0\psi\ .
\label{est:nlognneg}
\end{equation}
The previous computations \fer{est:nc}\,,\fer{est:nlognpos} and \fer{est:nlognneg} give us that the initial energy \eqref{eq:energyPE} is bounded as follows
$$
-C(m,\|e^{-\psi}\|_{L^1(\R^d)})-\int_{\R^d}n_0 \psi
-C(M,d)\|n_0\|_{L^\f d2(\R^d)}\le
\mathcal E[n_0]\le\|n_0\|_{L^\f d2(\R^d)}\,.
$$

Now, let us extend the above estimates locally in time. Estimate \fer{est:nlognpos} is obviously true for $n(t)$. Concerning the weighted $L^1$-norm of $n$ we have
$$
\dfrac{d}{dt}\int_{\R^d}n \psi =-\int_{\R^d} n\nabla \psi \cdot\nabla(\log n - c) 
\le \f1{4\delta} \int_{\R^d} \big|\nabla \psi \big|^2 n + 
\delta\int_{\R^d} n|\nabla (\log n - c)|^2\ ,
$$
with arbitrary $\delta>0$ to be choosen later. Hence, it follows that
\begin{equation}
\int_{\R^d}n(t) \psi \le \int_{\R^d}n_0 \psi 
+ \f1{4\delta}\|\nabla \psi \|_{L^\infty(\R^d)}^2 M\,t
+\delta\int_0^t\int_{\R^d} n|\nabla (\log n - c)|^2\ .
\label{est:momentum}
\end{equation}
It remains to estimate the energy dissipation from \eqref{eq:energy dissipation}. Using \fer{est:nc}, \fer{est:nlognneg} and \eqref{eq:energyPE}, we obtain
\begin{equation}
\begin{array}{rcl}
\int_0^t\int_{\R^d}\!\!\!\!\!\!\! &&n(s)|\nabla (\log n(s) - c(s))|^2\ dx\ ds
=\mathcal E[n_0]-\mathcal E[n](t)\\ \\
&&\le\|n_0\|_{L^\f d2(\R^d)}+
\int_{\R^d}(n(t)\log n(t))_-\, +\f12\int_{\R^d} n(t)c(t)\\ \\
&&\le\|n_0\|_{L^\f d2(\R^d)}
+C(m,\|e^{-\psi}\|_{L^1(\R^d)})+\int_{\R^d}n(t) \psi 
+C(M,d)\|n(t)\|_{L^{\f d2}(\R^d)}\ .
\end{array}
\label{est:dissipation}
\end{equation}
Plugging now \fer{est:momentum} with any $\delta\in(0,1)$ into \fer{est:dissipation}, we finally get
\begin{equation}
(1-\delta)\int_0^t\int_{\R^d}n |\nabla (\log n  - c )|^2 
\le C(n_0,\psi,M)(1+t)+C(M,d)\|n(t)\|_{L^{\f d2}(\R^d)}\ .
\label{est:dissipation2}
\end{equation}
Consequently, from \fer{est:momentum}, the same estimate holds for the weighted $L^1$-norm of $n$ , i.e.
\begin{equation}
\int_{\R^d}n(t) \psi \le 
C(n_0,\psi,M)(1+t)+C(M,d)\|n(t)\|_{L^{\f d2}(\R^d)}
\label{est:momentum2}
\end{equation}
and this gives us the control on the energy
\begin{equation}
-C(n_0,\psi,M)(1+t)-C(M,d)\|n(t)\|_{L^{\f d2}(\R^d)}\le
\mathcal E[n](t)\le \|n(t)\|_{L^{\f d2}(\R^d)}\,,
\label{est:energy}
\end{equation}
and exactly the same control on the entropy.\\

\begin{remark}[Logarithmic decay of the density]
The same estimates of this subsection hold true if we assume that $n_0\in L^1(\R^d,|x|^2dx)$ but with slightly modified computations (see \cite{BDP,CC08}). Using here the hypothesis $n_0\in L^1(\R^d, \psi(x)dx)$ we want to underline that one can obtain local and consequently global existence results under weaker condition on the decay at infinity of the initial density. For instance, it is sufficient to consider logarithmic decay at infinity for $n_0$, i.e. $\psi(x)=d\log(1+|x|^2)$.
\label{rm:logdecay}
\end{remark}

\subsection{Proof of Theorem \ref{the:localexistence}}
\label{prooflocalexistence}
To prove Theorem \ref{the:localexistence} we need to regularize system \fer{KS} with $\eps=0$ and $c$ given by \fer{eq:chemical} and to show that the sequence of regularized solutions satisfies firstly the previous estimates and secondly is relatively compact in an appropriate topological  space. This is now a quite standard procedure, used in the context of Keller-Segel system for example in \cite{BCM,CC08,CP08,CPZ04}. However, we give here a quite rapid sketch of the proof for the sake of completeness.

As a regularized problem we consider 
\begin{equation}
\partial_t\, n^\sigma=\Delta n^\sigma - \nabla\cdot(n^\sigma\nabla c\,^\sigma)\, ,
\label{KSregular1}
\end{equation}
where $c\,^\sigma$ is given by
\begin{equation}
c\,^\sigma(x,t)=\left\{
\begin{array}{rcl}
(E_d*n^\sigma(t)*\rho^\sigma)(x)\,,\quad\alpha=0\,,  \medskip\\
(B_d^\alpha*n^\sigma(t)*\rho^\sigma)(x)\,,\quad\alpha>0\,,
\end{array}
\right.
\label{KSregular2}
\end{equation}
and $\rho^\sigma$ is some sequence of smooth positive mollifiers with $\|\rho^\sigma\|_{L^1(\R^d)}=1$. The regularized initial condition
\begin{equation}
n_0^\sigma=n_0*\rho^\sigma
\label{ICregular}
\end{equation}
is also considered. Problem \fer{KSregular1}-\fer{KSregular2}-\fer{ICregular} has a nonnegative smooth solution as follows by the Schauder's fixed-point theorem. Moreover, the solution $(n^\sigma,c^\sigma)$ satisfies all the a priori estimates given in the previous sections. Indeed, one has essentially to check that the master equation \fer{est:np}  holds true at least as an inequality, so that the fundamental estimates \fer{est:n(p+1)} and \fer{est:n(p+1)bis} on the $L^{p+1}$ norm of $n^\sigma$ can be applied. This is the case since
$$
-\int_{\R^d}(n^\sigma)^p\,\Delta c^\sigma\le\int_{\R^d}(n^\sigma)^p\,(n^\sigma*\rho^\sigma)\le\int_{\R^d}(n^\sigma)^{p+1}\ .
$$

Concerning the compactness of the sequence $\{n^\sigma\}$, we are intended to use the Aubin compactness lemma in \cite{AU}. Therefore, we chose $B=L^2(\R^d)$, $X=H^1(\R^d)\cap L^2(\R^d, \psi(x)^{\f12}dx)$ compactly imbedded in $B$ and $Y=H^{-1}(\R^d)$ so that $B\subset Y$.  Using the previous a priori estimates, we prove below that $\{n^\sigma\}$ is bounded in $L^2((\delta,T_a-\delta);X)$ uniformly in $\sigma$, where $\delta\in(0,T_a)$ is arbitrarily small and $T_a$ is defined in \eqref{eq:Ta}. The boundedness of $\{\partial_t\,n^\sigma\}$ in $L^2((\delta,T_a-\delta);Y)$ uniformly in $\sigma$ follows  then from the same computation. Consequently, the Aubin's lemma gives that $\{n^\sigma\}$ is relatively compact in $L^2((\delta,T_a-\delta);B)$ and the proof is complete. Let us observe that for high dimensions ($d\ge6$) we obtain in fact that $\{n^\sigma\}$ is relatively compact in $L^2((0,T_a-\delta);B)$.

For the sake of simplicity, we omit the index $\sigma$ in the sequel. First of all, let us observe that $n\in L^\infty((\delta,T_a-\delta)\,;L^2(\R^d))$ if $d=3$ and $n\in L^\infty((0,T_a-\delta)\,;L^2(\R^d))$ if $d\ge4$, as it follows by \fer{est:hyper} and by \fer{est:nd/2+bis} respectively, being $a>\f d2$. 

Next, from \fer{est:gradc} with $p=\f d2$, we have that $|\nabla c|\in L^\infty((0,T_a-\delta)\,;L^d(\R^d))$, for any $d\ge3$. Consequently, since by  the H\"older's inequality, 
$$
\|(n\,\nabla c)(t)\|_{L^2(\R^d)}\le \|n(t)\|_{L^{\f{2d}{d-2}}(\R^d)}\|\nabla c(t)\|_{L^d(\R^d)}\,,
$$
using again \fer{est:hyper} and \fer{est:nd/2+bis} we obtain that
$n|\nabla c|\in L^\infty((\delta,T_a-\delta)\,;L^2(\R^d))$ for $d=3,4,5$ and $n|\nabla c|\in L^\infty((0,T_a-\delta)\,;L^2(\R^d))$ for $d\ge6$. 

Finally, multiplying the equation on $n$ in \fer{KS} against $n$, integrating over $\R^d$ and then over $(\delta,T_a-\delta)$, one easily obtain
$$
\int_\delta^{T_a-\delta}\|\nabla n(s)\|_{L^2(\R^d)}^2\, ds\,\le \|n(\delta)\|_{L^2(\R^d)}^2
+\int_\delta^{T_a-\delta}\|(n\,\nabla c)(s)\|_{L^2(\R^d)}^2\, ds\,,
$$
i.e. $|\nabla n|\in L^2((\delta,T_a-\delta)\,;L^2(\R^d))$ for $d=3,4,5$ and $|\nabla n|\in L^2((0,T_a-\delta)\,;L^2(\R^d))$ for $d\ge6$.

It remains to estimate the $L^2(\R^d, \psi(x)^{\f12}dx)$ norm of $n$. This is an immediate consequence of the computation
$$
\int_{\R^d} \psi ^{\f12}\,n^2(t)\le\left(\int_{\R^d} \psi \,n(t)\right)^{\f12}\,
\left(\int_{\R^d}n^3(t)\right)^{\f12}\ ,
$$
and of estimates \fer{est:momentum2} and again \fer{est:hyper} and \fer{est:nd/2+bis}. Therefore, $n\in L^\infty((\delta,T_a-\delta)\,;L^2(\R^d, \psi(x)^{\f12}dx))$ for $d=3,4,5$ and $n\in L^\infty((0,T_a-\delta)\,;L^2(\R^d,\psi(x)^{\f12}dx))$ for $d\ge6$. 

\hfill\cqfd

\begin{remark} [Improved smallness condition]
\label{bestsmallnesscond}
In {\rm \cite{CPZ04}}, in order to obtain the global existence result, the authors have used  inequality \eqref{est:n(p+1)Sob} instead of \eqref{est:n(p+1)bis} and the corresponding inequality
\begin{equation}
\f{d}{dt}\int_{\R^d} n^p(t)\le(p-1)\Big[C_{S}^2(d) \norm{n(t)}_{L^{\frac d 2}(\R^d)}-\f4p\Big]\norm{\nabla n^{\f p2}}_{L^{2}(\R^d)}^2\ .
\label{est:nd/2bis}
\end{equation}
instead of \eqref{est:nd/2}, with the corresponding smallness condition  
\begin{equation}
\norm{n_0}_{L^{\frac d 2}(\R^d)}<\,\f{8}{d\,C^2_S(d)}\ .
\label{n_0condSob}
\end{equation}
Since we can prove that $C_{GN}^{\,2(1+\f2d)}\left(\f d2\,,d\right)< C_S^2(d)$, condition \eqref{n_0condGN} is a weaker condition than \fer{n_0condSob}. The last together with the facts that
\begin{enumerate}[(i)]
\item inequality \fer{est:n(p+1)bis} is sharp (equality holds true for the minimizers of \fer{GN2})\,;
\item under the smallness condition \fer{n_0condSob} : all the $\norm{n(t)}_{L^{p}(\R^d)}$ norms are time decreasing,  for $p\in(1,\min\{p^*,d\}]$ and where $p^*:=4\,\norm{n_0}^{-1}_{L^{\frac d 2}(\R^d)}C_{S}^{-2}(d)$ is greater than $\f d2$; the hyper-contractivity property \eqref{est:hyper} holds true globally in time for $p>\min\{p^*,d\}$; the entropy $\int_{\R^d} n\log n$ is time decreasing in dimension $d=3$;
\end{enumerate}
induced to conjecture that the smallness condition \fer{n_0condGN} is the optimal one.

To prove the claimed inequality between optimal constants, it is sufficient to write the Gagliardo-Nirenberg inequality \fer{GN2} for $\sigma=\f2d$ and the function $v$ for which the equality holds true, i.e.
$$
C_{\f2d,\,d}^{\,2(1+\f2d)}\|\nabla v\|_{L^2(\R^d)}^2\,\|v\|_{L^2(\R^d)}^{\f4d}=\|v\|_{L^{2(1+\f2d)}(\R^d)}^{2(1+\f2d)}\ .
$$
Observing that by standard interpolation and the Sobolev inequality \fer{in:Sob} we have
$$
\|v\|_{L^{2(1+\f2d)}(\R^d)}^{2(1+\f2d)}\le
\|v\|_{L^2(\R^d)}^{\f4d}\,\|v\|_{L^{\f{2d}{d-2}}(\R^d)}^2
\le C_S^2(d)\,\|\nabla v\|_{L^2(\R^d)}^2\,\|v\|_{L^2(\R^d)}^{\f4d}\ ,
$$
and using \fer{GNcost}, the inequality $C_{GN}^{\,2(1+\f2d)}\left(\f d2\,,d\right)\le C_S^2(d)$ is proved. Since equality cannot hold true in \eqref{est:n(p+1)Sob} while it can in \eqref{est:n(p+1)bis}, the above claimed strict inequality is obtained. 

The time decreasing behaviour of the family of $L^p$-norms of $n$, $p\in(1,\min\{p^*,d\}]$, whenever the smallness condition \fer{n_0condSob} holds true, is a direct consequence of inequality \eqref{est:nd/2bis}, where the constant $C_{S}^2(d)$ doesn't depend on $p$. In order to obtain the same result under the improved smallness condition \eqref{n_0condGN}, one should know the behaviour of the constant $C_{GN}^{\,\f 2p(p+1)}(p,d)$ in \eqref{est:nd/2} with respect to $p$. 

The hyper-contractivity property \eqref{est:hyper} holds true globally in time for $p>\min\{p^*,d\}$ since we can choose the modulus of ``$\f d2$-equintegrability'' in \eqref{equi} time independent. 

Finally, let us show that \fer{n_0condSob} implies the time decreasing behaviour of the entropy $\int_{\R^d} n\log n$, at least in dimension $d=3$. Indeed, 
\begin{equation}
\f{d}{dt}\int_{\R^d} n\log n=-4\int_{\R^d}|\nabla\sqrt n|^2 +\int_{\R^d} n^2 \ .
\label{est:nlogn}
\end{equation}
Then, using the Gagliardo-Nirenberg inequality \fer{GN} with $p=1$ and $v=\sqrt n$ to estimate $\int_{\R^d} n^2$ in \fer{est:nlogn}, we have
\begin{equation}
\f{d}{dt}\int_{\R^d} n\log n\le\norm{\nabla\sqrt n}^2_{L^2(\R^d)}\Big[C^4_{GN}(1,d)\norm{n}_{L^{\f d2}(\R^d)}-4\Big]\ .
\label{est:nlogn2}
\end{equation}
The best constant $C_{GN}(1,d)$ is unknown (at our best knowledge) except for $d=3$ (see \cite{DelPinoDolbeault}). In that case $C_{GN}(1,3)=(\f1{2\pi^2})^{1/6}$, while $C_S^2(3)=\f43(\f1{2\pi^2})^{2/3}$ from \fer{in:Sob} and the issue is proved since the right hand side of \fer{est:nlogn2} turns to be negative under the smallness condition \fer{n_0condSob}. 
\end{remark}
\begin{remark}
\label{paraboliccase}
In the parabolic-parabolic case ($\eps>0$) we can not have the local existence result without smallness condition on the initial densities, obtained for the parabolic-elliptic system. This is due essentially to the fact that in the master equation \fer{est:np} one has to add the term $-\eps\,(p-1)\int_{\R^d}n^p\,\partial_tc$ stronger with respect to $\int_{\R^d}n^{p+1}$, since the only control that we have for the temporal derivative of $c$ is an $L^2$ control from \eqref{eq:energy dissipation}. For local and global existence results of \fer{KS} with $\eps>0$, we refer to~\cite{CP08, KozoSugy08} and the references therein.
\end{remark}
%


\section{Blow-up for the parabolic-elliptic system ($\eps=0$)}
\label{Sec:Blow-up}
In order to characterize blowing up solutions of the parabolic-elliptic \fer{KS} system, the general idea is to follow the evolution of the second moment of $n$, i.e. $I(t):=\int_{\R^d}|x|^2n(x,t)\ dx$, and to prove that $I$ satisfies a differential inequality of type
\begin{equation}
\f d{dt}I(t)\le f(I(t))\ ,
\label{eq:I'(t)}
\end{equation}
where $f$ is a  continuous nondecreasing function such that $f(0)<0$. Indeed, whenever we can exhibit \eqref{eq:I'(t)}, defining $I^*:=\inf\{I>0\,|\,f(I)=0\}\in (0,+\infty]$, for any sufficiently smooth solution $n$ of \fer{KS} with finite initial second moment satisfying $I(0)<I^*$, there exists a time $T^*<\infty$ such that $\lim\limits_{t\nearrow T^*}I(t)=0$. The vanishing behaviour of the second moment implies the blow-up of some norm critical for the existence of the solution and henceforth of the solution itself.  Obviously the blowing time is in general smaller than $T^*$ and at our best knowledge it is an open problem to characterize the blowing time in term of the exploding critical norm instead of $I(0)$.

The above technique as been applied firstly by Biler in \cite{B95} for a model of gravitational interaction of particles on a star-shaped domain of $\R^d$, similar to the \fer{KS} system with $\eps=\alpha=0$. Successively, it has been used by several authors in the context of the Keller-Segel system (see for instance \cite{Nagai95,CPZ04,BDP}). The methodology is also reminiscent of the blow-up criteria for the nonlinear Schr\"{o}dinger equation initiated by Glassey \cite{G}, (see also \cite{W}), and has been successively applied to kinetic gravitational models \cite{GlasseyBook} and kinetic chemotaxis models \cite{BournaveasC}.

Notice that there exists an alternative (and non constructive) method to obtain the existence of blowing up solutions, based on energy features (unbounded from below) and the analysis of the possible stationary states. We refer to \cite{Horstmann02,HW05} for more details. 

Concerning the parabolic-elliptic \fer{KS} system in two space dimensions, some candidate function $f$ in \fer{eq:I'(t)} can be explicitly computed, whatever $\alpha\ge0$ is, and the condition $f(0)<0$ reads as $M>8\pi$ when  $\alpha = 0$, (see \cite{BDP,CC08}). This threshold for the mass $M$ is sharp since global existence can be proved for $M<8\pi$ and suitable conditions on the initial density $n_0$, (\cite{BDP, CC08}).  When blow-up occurs, the solution ceases to exist classically and possible extensions after the blow-up time have been proposed in \cite{VelazquezA,VelazquezB,DolbeaultSchmeiser}, depending upon the system regularization. Finally, in the critical case $M=8\pi$, it has been proved in \cite{BCM} that weak free energy solutions still exist globally in time and that they blow up in infinite time.

In the case $d\ge3$ the derivation of \fer{eq:I'(t)} is a little more complicated since the potential appears in the evolution equation for the moment $I$. More specifically, from 
\begin{equation}
\dfrac{d }{dt}I(t)= 2d M+2\int_{\R^d}n(x,t)\,x\cdot\nabla c(x,t)\, dx\,,
\label{eq:evol I} 
\end{equation}
we have (after symmetrization of the integral term) for $\alpha=0$
\begin{equation}
\dfrac{d }{dt}I(t)=2dM-\dfrac1{|\Sph^{d-1}|}\iint_{\R^d\times\R^d}n(x,t)\f{1}{|x-y|^{d-2}}n(y,t)\,dxdy
\label{eq:I1} 
\end{equation}
while for $\alpha>0$
\begin{equation}
\dfrac{d }{dt}I(t)=2dM-\dfrac1{|\Sph^{d-1}|}\iint_{\R^d\times\R^d}n(x,t)\f{g_\alpha(|x-y|)}{|x-y|^{d-2}}n(y,t)\,dxdy\,,
\label{eq:I2}
\end{equation}
where $g_\alpha$ is defined in Lemma \ref{lem:Besselproperties}. Therefore, different blow-up criteria can be obtained according to how the right hand side of \eqref{eq:I1} or \eqref{eq:I2} is estimated with respect to $I$, $\mathcal E[n]$ and $M$. Let us observe that \fer{eq:I2} becomes \fer{eq:I1} when $\alpha=0$ since $g_0(|x|)\equiv1$. The continuity with respect to the dimension $d$ is also satisfied in the sense that from \fer{eq:I1} and \fer{eq:I2} we obtain the evolution equation of $I$ when~$d=2$.

In the sequel we shall derive two inequalities of type \eqref{eq:I'(t)} together with the consequent blow-up criteria  and we shall discuss their complementarity. One of the criteria is yet contained in \cite{CPZ04} and both are contained in \cite{B95} for the case of star-shaped domains. However, only the case of the parabolic-elliptic  \eqref{KS} system in absence of chemical degradation ($\eps=\alpha=0$) was considered in the cited papers. Here we improve the criterion involving the free energy (not required to be negative) and we extend both to the case $\alpha>0$. 

\subsection{Derivation of two blow-up criteria}
\begin{proposition}[First  blow-up criterion]
Let $\eps=0$, $\alpha\ge0$, $d\ge3$ and $a>\f d2$. Assume that the nonnegative initial density $n_0\in (L^1\cap L^a)(\R^d)$ has finite second momentum satisfying
\begin{equation}
\int_{\R^d} |x|^2 n_0(x)\, dx < K_1^\alpha(d,M)  M^{\frac{d}{d-2}}\,, 
\label{1BU}
\end{equation}
where  $K_1^\alpha(d,M)$ is defined in \eqref{eq:K1alpha} for $\alpha>0$ and 
\begin{equation}
K_1^0(d,M)= K_1(d):=2^{-\f d{d-2}}(d |\Sph^{d-1}|)^{-\f2{d-2}}\,.
\end{equation}
Then, the solution to the parabolic-elliptic \eqref{KS} system constructed in Theorem \ref{the:localexistence} blows up in finite time,  that is the maximal time of existence $T_{max}$ is finite and $\lim\limits_{t \nearrow T_{max}}\|n(t)\|_{L^a(\R^d)}=+\infty$. 
Moreover, the blow-up condition \fer{1BU} is incompatible with the smallness condition \fer{n_0condSob} for global existence. 
\label{prop:1BUcriterion}
\end{proposition}
\begin{proposition}[Second blow-up criterion]
Let $\eps=0$, $\alpha\ge0$, $d\ge3$ and $a>\f d2$. Assume that the nonnegative initial density $n_0\in (L^1\cap L^a)(\R^d)$ has finite second momentum satisfying
\begin{equation}
\int_{\R^d} |x|^2 n_0(x)\, dx < K_2(d)\, M^{\,1+\frac{2}d}\, \exp\left(-\dfrac2{dM}\,\mathcal{E}[n_0]\right)\,,  \qquad K_2(d):=\f d{2\pi}\,e^{-\f d{d-2}}\,.
\label{2BU}
\end{equation}
Then, the solution to the parabolic-elliptic \eqref{KS} system constructed in Theorem \ref{the:localexistence} blows up in finite time, that is the maximal time of existence $T_{max}$ is finite and $\lim\limits_{t \nearrow T_{max}}\|n(t)\|_{L^a(\R^d)}=+\infty$. Moreover, the blow-up condition \fer{2BU} are incompatible with the smallness condition \eqref{n_0condSob} for global existence. 
\label{prop:2BUcriterion}
\end{proposition}
\begin{proof}[Proof of Proposition \ref{prop:1BUcriterion}]
Assume first for simplicity that $\alpha=0$. Using the left inequality in \fer{est:potential} into \fer{eq:I1}, we have:
\[
\dfrac{d }{dt}I(t)\le 2dM-|\Sph^{d-1}|^{-1}2^{\,1-\f d2}\,M^{\,\f d2+1}\,I^{\,1-\f d2}(t)\,,
\]
i.e. \fer{eq:I'(t)} holds true with the increasing function $f(\lambda)=2dM-|\Sph^{d-1}|^{-1}2^{\,1-\f d2}\,M^{\,\f d2+1}\,{\lambda}^{\,1-\f d2}$ satisfying $f(\lambda)\to-\infty$ as $\lambda\searrow0$ and $f(I^*)=0$ with $I^*:=K_1(d)\, M^{\frac{d}{d-2}}$. Hence, we obtain an obstruction to global existence when $I(0)<I^*$, which is the condition \eqref{1BU}. From estimates \eqref{est:entropy2} and \eqref{est:potential}, we deduce the existence of $T_1>0$, $T_2>0$ and $T_{max}>$ such that
\[
\lim\limits_{t\nearrow T_1}\int_{\R^d}n\log n=+\infty\,,
\quad
\lim\limits_{t\nearrow T_2}\int_{\R^d}nc=+\infty
\quad\hbox{and}\quad
\lim\limits_{t\nearrow T_{max}}\|n(t)\|_{L^a(\R^d)}=+\infty\,
\]
with $T_{max}\le\min\{T_1,T_2\}$.

When $\alpha>0$, in order to use the potential confinement Lemma as before into \eqref {eq:I2}, because of the  decreasing behaviour of $g_\alpha$, one has to proceed as following
\begin{align}
\dfrac d{dt}I(t)&=2dM-\f1{|\Sph^{d-1}|}\iint_{\R^d\times\R^d}n(x,t)\f{g_\alpha(|x-y|)}{|x-y|^{d-2}}n(y,t)\,dxdy
\nonumber\\
&\le 2dM-\f{g_\alpha(R)}{|\Sph^{d-1}|}\iint_{|x-y|<R}n(x,t)\f1{|x-y|^{d-2}}n(y,t)\,dxdy \nonumber\\
&\le 2dM-\f{g_\alpha(R)}{|\Sph^{d-1}|}\iint_{\R^d\times\R^d}n(x,t)\f1{|x-y|^{d-2}}n(y,t)\,dxdy
\nonumber\\
&\quad\quad+\f{g_\alpha(R)}{|\Sph^{d-1}|} \iint_{|x-y|> R}n(x,t)\f1{|x-y|^{d-2}}n(y,t)\,dxdy\,.
\nonumber
\end{align}
Then, using the left inequality in \fer{est:potential}, we obtain for any $R>0$
\begin{equation}
\dfrac d{dt}I(t)\le 2dM-\f{g_\alpha(R)}{|\Sph^{d-1}|}\,2^{1-\f d2} M^{\,\f d2+1} I^{\,1-\f d2}(t)+\f{g_\alpha(R)}{|\Sph^{d-1}|\,R^{d-2}} M^2\,.
\label{eq:I'new}
\end{equation}
The differential inequality \eqref{eq:I'new} is again of type \eqref{eq:I'(t)}, with in the r.h.s. an increasing function of $I$, converging to $-\infty$ as $I \searrow0$ and to a positive constant as  $I$ goes to $+\infty$. To obtain the best   blow-up condition, one has to maximize with respect to $R$ the value $I^*(R)$ for which the r.h.s. of \eqref{eq:I'new} vanishes.  Such a computation give the definition $I^*:=K_1^\alpha(d,M)  M^{\frac{d}{d-2}}$, where
\begin{equation}
K_1^\alpha(d,M):=\f12\left[\sup_{R>0}\left(\f{g_\alpha(R)\,R^{d-2}}{2d|\Sph^{d-1}|\,R^{d-2}+g_\alpha(R)\,M}\right)
\right]^{\f2{d-2}}\,.
\label{eq:K1alpha}
\end{equation}
Owing to the properties of $g_\alpha$ for $\alpha>0$, the function to maximise into formula \eqref{eq:K1alpha} is a positive continuous function, vanishing for $R=0$ and decaying to $0$ as $R\to\infty$. Therefore, the supremum in \eqref{eq:K1alpha} is achieved and $K_1^\alpha$ is well defined. Consequently, the blow-up criterion \eqref{1BU} has been obtained. Moreover, $K_1^0(d,M) = K_1(d)$ because $g_0(|x|)\equiv 1$ and the supremum in \eqref{eq:K1alpha} is achieved for $R\to+\infty$.

Finally, let us prove the claimed incompatibility. It is sufficient to consider the case $\alpha=0$ since $K_1^\alpha$ is strictly decreasing with respect to $\alpha$ so that $K_1^\alpha(d,M)<K_1(d)$. From the potential confinement Lemma \ref{lem:estpotential} and from \fer{1BU}, we obtain
\begin{equation}
\|n_0\|_{L^{\f d2}(\R^d)}\ge 2^{1-\f d2}\,M^{\f d2}\,C_{HLS}^{-1}(d,d-2)\,I^{1-\f d2}(0)>
2\,d\,|\Sph^{d-1}|\,C_{HLS}^{-1}(d,d-2)\ .
\label{estincomp1}
\end{equation}
Using the relation \fer{eq:CS-HLS} into \fer{estincomp1}, the reverse of condition \fer{n_0condSob} follows. 

\end{proof}
\begin{proof}[Proof of Proposition \ref{prop:2BUcriterion}]
Let us again consider firstly the case $\alpha=0$. Then, we shall make use of the definition of the free energy, its time decreasing behaviour and of  the entropy lower bound \fer{est:entropy2} instead of the potential confinement inequality into \eqref{eq:I1},  to get
\begin{align}
\dfrac{d }{dt}I(t)& = 2dM-(d-2)\int_{\R^d} n(x,t)c(x,t)\, dx\nonumber \\
& =
2dM+2(d-2) \left(\mathcal E[n](t) -\int_{\R^d} n(x,t)\log n(x,t)\, dx \right) 
\label{eq:BUargument}\\ 
& \leq d(d-2) M\log I(t) +2(d-2)  \mathcal E[n_0] + B(d,M)
\, , \nonumber
\end{align}
where the constant $B(d,M)$ is defined as
\begin{equation}
B(d,M) := d^2M-2(d-2)\left[M\log M+\f{dM}2\log\left(\f{dM}{2\pi}\right)\right]\, .
\label{eq:bdm}
\end{equation}
Therefore, \fer{eq:I'(t)} yields with the nondecreasing function $f(\lambda)= d(d-2) M\log \lambda+2(d-2)  \mathcal E[n_0] + B(d,M)$, satisfying $f(\lambda)\to-\infty$ as $\lambda \searrow0$ and $f(I^*)=0$ with $I^*:=\exp\left(-\f 2{dM} \mathcal E[n_0] - \f{B(d,M)}{d(d-2) M} \right)$. Again, for $I(0)<I^*$, that is condition \eqref{2BU}, we obtain an obstruction to global existence as before.

When $\alpha>0$, we proceed as for $\alpha=0$, but using formula \eqref{eq:Euler alpha} for the Bessel kernel into the evolution equation of $I$ write as \eqref{eq:evol I}. Therefore, after the symmetrization of the integral term in \eqref{eq:evol I}, we have
\begin{eqnarray*}
\dfrac{d }{dt}I(t) &=&  2dM+\iint_{\mathbb{R}^{d}\times\mathbb{R}^{d}} n(x,t)(x-y)\cdot\nabla B_d^\alpha(x-y)\,n(y,t)\, dxdy\\
&= & 2dM-(d-2)\iint_{\mathbb{R}^d\times\mathbb{R}^d} n(x,t)\,B_d^\alpha(x-y)\,n(y,t)\, dxdy \nonumber\\
 &~&  -2\alpha\iint_{\mathbb{R}^d\times\mathbb{R}^d} n(x,t)(B_d^\alpha*B_d^\alpha)(x-y)\,n(y,t)\, dxdy\\
 &=& 2dM-(d-2)\int_{\mathbb{R}^d} n(x,t)c(x,t)dx-2\alpha\int_{\mathbb{R}^d}(B_d^\alpha*n)^2(x,t)dx\\
 &=& 2dM+2(d-2)\left(\mathcal E[n](t)-\int_{\R^d} n(x,t)\log n(x,t)\,dx \right)-2\alpha\int_{\mathbb{R}^d}c^2(x,t)\, dx\, .
\end{eqnarray*}
Neglecting the last negative contribution we are reduced to the previous estimate \eqref{eq:BUargument} and we can conclude as above.

Finally, let us prove the claimed incompatibility. Applying the logarithmic function to both sides of  \fer{2BU} and using the definition of the energy, we arrive easily to the inequality:
\begin{equation}
\f{dM}2\left(\log I(0)+1\right)+\int_{\R^d}n_0(x)\log n_0(x)\, dx-M\log M-\f{dM}2\log\left(\f{dM}{2\pi}\right)
<\f12\int_{\R^d}n_0(x)c_0(x)\, dx-\f{dM}{d-2}\ .
\label{estincomp3}
\end{equation}
The left hand side of \fer{estincomp3} is nonnegative owing to Lemma \ref{est:entropy}. Then, from \fer{estincomp3} and again the potential confinement Lemma \ref{lem:estpotential}, we have
\begin{equation}
\f{2dM}{d-2}<\int_{\R^d}n_0(x)c_0(x)\, dx\le\mu_d\,C_{HLS}(d,d-2)\,M\,\|n_0\|_{L^{\f d2}(\R^d)}\,,
\label{estincomp2}
\end{equation}
and the reverse of condition \fer{n_0condSob} follows.

\end{proof}
\begin{remark}
We are actually not able to prove that both criteria, \eqref{1BU} and \eqref{2BU}, are also incompatible with the improved smallness condition \fer{n_0condGN} for the global existence. This is due to the lack of knowledge about the constant $C_{GN}(d/2,d)$.
\end{remark}
\begin{remark}
(i) It is easy to check the invariance of criterion \eqref{1BU} for $\alpha=0$ and of criterion  \eqref{2BU} for any $\alpha\ge0$, under the scaling $n_0(x)\to n_0^\lambda(x)=\lambda^{-2}n_0(\lambda^{-1}x)$\,, preserving the $L^{\f d2}-$norm. This invariance doesn't hold true for criterion \eqref{1BU} when $\alpha>0$ because of the dependency on $M$ of the constant $K_1^\alpha$. However, since $K_1^\alpha(d,M)$ is strictly decreasing with $\alpha$\,, so that $K_1^\alpha(d,M)<K_1(d)$ for any $\alpha>0$, and because of the claimed incompatibility, it is not possible to construct from a density $n_0$ satisfying the smallness condition \eqref{n_0condSob} a density $n_0^\lambda(x)=\lambda^{-2}n_0(\lambda^{-1}x)$ satisfying criterion \eqref{1BU} for some $\alpha>0$. (ii) The decreasing behaviour of $K_1^\alpha$ with respect to $\alpha$ is naturally due to the fact that the degradation term $-\alpha c$ in the \eqref{KS} system prevents blow-up of the density $n$. This term has been completely neglected in the deduction of the second blow-up criterion \eqref{2BU} involving the energy. It would be interesting, eventhough  hard to do, not to neglect it.
\end{remark}
\begin{remark}
Looking attentively at the proofs of Propositions \ref{prop:1BUcriterion} and  \ref{prop:2BUcriterion}, it appears evident that the second blow-up criterion is sharper than the first one. Indeed, the first blow-up criterion has been obtained with the help of the potential confinement inequality \eqref{est:potential}, while the second blow-up criterion has been obtained applying the entropy lower bound \eqref{est:entropy2} and that the latter is a sharp inequality (see Remark \ref{rem:sharpentropylowerbound}).
\label{rk:bestcriterion}
\end{remark}
\begin{remark}[Corrected energy]
Coming back to the line \eqref{eq:BUargument} and using the fact that $\f d{dt}\mathcal E[n(t)]\le0$ we deduce the following differential inequality for the corrected energy defined in \eqref{eq:corrected energy}
\begin{equation}
I(t)\f d{dt}\mathcal F[n](t) \le d(d-2) M \mathcal F[n](t)+B(d,M)\,,
\label{eq:correctedenergyinequality}
\end{equation}
Moreover, the blow-up criterion involving the initial free energy \eqref{2BU} reads equivalently as: 
\begin{equation}
d(d-2) M \mathcal F[n_0]+B(d,M)<0\,.
\label{eq:2BUequiv}
\end{equation}
In Section \ref{Sec:Blow-up para-para} we will generalize inequality  \eqref{eq:correctedenergyinequality} as well as the blow-up condition \eqref{eq:2BUequiv} to the parabolic-parabolic Keller-Segel system \eqref{KS} in order to obtain a concentration result for the $L^{\f d2}$-norm of $n$. 
\label{rk:correctedenergy}
\end{remark}
\subsection{Complementarity of the blow-up criteria}
\label{sec:complementarity}
The goal of this subsection is to construct examples of initial data $n_0$ satisfying  either
\begin{equation}
K_2(d)\, M^{\,1+\frac{2}d}\, \exp\left(-\dfrac2{dM}\,\mathcal{E}[n_0]\right) < \int_{\R^d} |x|^2 n_0(x)\, dx < K_1(d)\, M^{\frac{d}{d-2}}\, ,
\label{eq:complement1}
\end{equation}
or
\begin{equation}
K_1(d)\, M^{\frac{d}{d-2}} < \int_{\R^d} |x|^2 n_0(x)\, dx < K_2(d)\, M^{\,1+\frac{2}d}\, \exp\left(-\dfrac2{dM}\,\mathcal{E}[n_0]\right)  \, .
\label{eq:complement2}
\end{equation}
Here $\alpha=0$. The case $\alpha>0$ is not considered for sake of simplicity.

We start with the latter \eqref{eq:complement2} which is more natural, since it is a direct consequence of the fact that the energy is unbounded from below and of the mixed homogeneities of the potential and entropy terms composing the energy with respect to the mass-preserving dilation $f(x)\to f^\lambda(x)=\lambda^{-d}f(\lambda^{-1}x)$. For this purpose, we introduce the following family of densities indexed by $\lambda>0$
\begin{equation}
n_0^\lambda(x) := \dfrac12\left[\lambda^{-d} \varphi\left(\dfrac {x-a}\lambda\right) + \lambda^{-d} \varphi\left(\dfrac {x+a}\lambda\right)\right]\, , 
\label{eq:density1}
\end{equation}
where $a\neq 0$ is some point to be chosen later and  $\varphi$ is a nonnegative function in $(L^1\cap L^{\f d2})(\R^d)$ such that $\int_{\R^d}\varphi(z)\, dz=M$  and $\mathrm{Supp}\,\varphi \subset B(0,1)$. Then, the densities $n_0^\lambda$ belong to $(L^1\cap L^{\f d2})(\R^d)$, have mass equal to $M$, $L^{\f d2}$-norm increasing as $\lambda\searrow0$ and the second moment given by
\begin{equation} 
\int_{\R^d}|x|^2n_0^\lambda(x)\,dx
= \dfrac12\left[\int_{\R^d} |a + \lambda z|^2\varphi( z)\, dz +  \int_{\R^d} |-a + \lambda z|^2\varphi( z)\, dz \right]= M |a|^2 + \lambda^2 \int_{\R^d} | z|^2\varphi( z)\, dz\,. 
\label{eq:equiv 1}
\end{equation}

When evaluating the free energy, the cross-interaction between the two densities located around $a$ and $-a$ is zero in the entropy term, if $\lambda$ is small enough so  that the supports $B(a,\lambda)$ and $B(-a,\lambda)$ are disjoints. Hence, we have
\begin{align}
\mathcal E[n_0^\lambda] & = \int_{\R^d} \varphi(z)\log \varphi(z)\, dz - dM \log \lambda-M\log2 - \lambda^{2-d}\ \f{\mu_d}4\iint_{\R^d\times\R^d} \varphi(z)\dfrac1{|z - z'|^{d-2}} \varphi(z')\, dzdz' \nonumber\\
& \qquad -\ \f{\mu_d}4 \iint_{\R^d\times\R^d} \varphi(z)\dfrac1{|2a + \lambda(z - z')|^{d-2}} \varphi(z')\, dzdz'\, ,
\label{eq:equiv 2}
\end{align}
 which goes to $-\infty$ as $\lambda\searrow0$. Comparing \eqref{eq:equiv 1} and \eqref{eq:equiv 2} clearly there exists $a\neq 0$ and $\lambda>0$ sufficiently small such that the corresponding density  $n_0^\lambda$ satisfies \eqref{eq:complement2}. 

 Notice that, the greater is $M$, the greater $|a|$ has to be chosen for the left inequality in \eqref{eq:complement2} to be satisfied. Nevertheless, in case of large mass $M$, if the two densities are concentrated enough, i.e. if  $\lambda$ is small enough, the potential in \eqref{eq:equiv 2} is strong enough to ensure  blow-up of the solution.

In order to exhibit an example of density satisfying \eqref{eq:complement1}, we follow the same idea as before but with the aim to obtain a corresponding energy $\mathcal E[n_0^\lambda] $ with the entropy term dominating the potential term. Therefore, we consider the following sequence of densities
\begin{equation}
n_0^\lambda(x) := \dfrac1N \sum_{i=1}^N \lambda^{-d}\varphi\left(\dfrac{x-a_i}\lambda\right)\,, \qquad \lambda = N^{1/(2 - d)}\, ,
\label{newsequence}
\end{equation}
where $\varphi$ is defined as above and in addition the family of points $(a_i)_{1\leq i\leq N}$ is choosen so that $a_i$ and $-a_i$ belong both to the family. Again, the densities $n_0^\lambda$ belong to $(L^1\cap L^{\f d2})(\R^d)$, have mass equal to $M$, $L^{\f d2}$-norm increasing as $\lambda\searrow0$ and the second moment given by
\begin{equation}
\int_{\R^d} |x|^2 n_0^\lambda(x)\, dx=
M\,\dfrac1N \sum_{i=1}^N|a_i|^2 + \lambda^2 \int_{\R^d}|z|^2\varphi(z)\, dz\, . 
\label{eq:mom}
\end{equation} 
Again we assume that $\lambda$ is chosen such that the supports $B(a_i,\lambda)$ of each contribution in \eqref{newsequence} are disjoints and we introduce the notation $D^\lambda(i,j)=\hbox{dist}(B(a_i,\lambda),B(a_j,\lambda))$. 
Computing separately each contribution of the energy functional, we obtain
\begin{equation}
\int_{\R^d} n_0^\lambda(x)\log n_0^\lambda(x)\, dx
= \dfrac1N \sum_{i=1}^N  \int_{\R^d} \varphi(z)\log \left(\dfrac{\lambda^{-d}}N \varphi(z)\right)\, dz=
-M\log(N\lambda^d)+\int_{\R^d} \varphi(z)\log\varphi(z)\ dz\,,
\label{eq:ent}
\end{equation} 
and
\begin{align}
& \iint_{\R^d\times\R^d} n_0^\lambda(x)  \dfrac1{|x-y|^{d-2}}  n_0^\lambda(y)\, dx dy \nonumber \\ 
& = \iint_{\R^d\times\R^d}\varphi(z) \f1{|z-z'|^{d-2} }  \varphi(z')\, dzdz' + \dfrac1{N^2}\sum_{i\neq j}\iint_{\R^d\times\R^d}\varphi(z) \f1{|a_i - a_j + \lambda(z-z')|^{d-2} }  \varphi(z')\, dzdz' \nonumber
\\
&\le \iint_{\R^d\times\R^d}\varphi(z) \f1{|z-z'|^{d-2} }  \varphi(z')\, dzdz'  + \dfrac{M^2}{N^2}\sum_{i\neq j} (D^\lambda(i,j))^{2-d} \,,
\label{eq:pot}
\end{align}
due to the choice $N\lambda^{d-2}=1$. 

We claim that there exists a family of points $(a_i)_{1\leq i\leq N}$ such that the last contribution in the r.h.s of \eqref{eq:pot} is uniformly bounded with respect to $N$. The argumentation goes as follows. First we may change the reference norm, i.e. we can replace the euclidean norm $|\cdot|_2$ in $\R^d$ with the supremum norm $|\cdot|_\infty$, up to some constant. Then, we distribute $N$ points on the regular grid $N^{-1/d}\cdot\Z^d$ inside the hypercube $[-1/2,1/2]^d$. Next, we observe that for any index $i$ and any integer $k<\f12 N^{1/d}$ there are at most $C k^{d-1}$ indices $j$ such that $|a_i - a_j|_\infty = N^{-1/d}  k$, where $C$ is a constant depending only on the dimension. As a matter of fact, after rescaling space by a factor $N^{1/d}$, those points are regularly distributed on a sphere of radius $k$. Finally, we have the following estimates as $N\to +\infty$: $\lambda  \ll N^{-1/d} $ and consequently $D^\lambda(i,j) \sim |a_i - a_j|_2$. To conclude, we can estimate the last contribution in the r.h.s of \eqref{eq:pot} as follow
\begin{equation*}
\dfrac1{N^2}\sum_{i\neq j} (D^\lambda(i,j))^{2-d} 
\lesssim C \dfrac1{N^2}\, N \sum_{k = 1}^{[2^{-1}N^{1/d}]} k^{d-1} \left( N^{-1/d} k \right)^{2-d}  \leq C \dfrac 1N\, N^{(d-2)/d}\,\f14\, N^{2/d} =\f 14\,C\,,
\end{equation*}
where $[2^{-1}N^{1/d}]$ denotes the integer part of $2^{-1}N^{1/d}$. Therefore the interaction potential is bounded from below, while the entropy \eqref{eq:ent} is decreasing towards $-\infty$ as $N\to +\infty$. Moreover, it is always possible to scale appropriately the location of the family (here inside $[-1/2,1/2]^d$ for the sake of reference) to ensure that the right inequality in \eqref{eq:complement1} is satisfied.

\begin{remark}[Exotic construction]
The above example \eqref{newsequence}, which is made of a large superposition of approximations of the idenity with disjoint supports, is not likely to be a configuration encountered genuinely. We will see later in Section \ref{Sec:num} an (over-)simplified version of the parabolic-elliptic \eqref{KS} system for which the second blow-up criterion \eqref{2BU} always contains the first one \eqref{1BU}.  
\end{remark}

\section{Concentration phenomenon for the parabolic-parabolic system}
\label{Sec:Blow-up para-para}
In this section we shall prove Theorem \ref{th:concentration}, i.e. we shall show that the parabolic-parabolic Keller-Segel system ($\varepsilon>0$) exhibits a concentration behaviour under some smallness condition on the  second moment of $n_0$, generalization of condition \eqref{2BU} obtained  in the parabolic-elliptic case. However, although we are able to obtain blow-up of the cell density in the latter case, we fail to do so as soon as $\e >0$, and we derive a weaker result claiming that the $L^{d/2}-$norm of cell density reaches high levels (at least for small $\varepsilon$).

Pretty few results are available for the parabolic-parabolic system \eqref{KS}. We stress out a recent contribution by Ci\'{e}slak and Lauren\c{c}ot \cite{CL}  in which they prove the occurrence of such parabolic-parabolic blow-up events in one space dimension but with suitable nonlinear diffusion. They follow a perturbation strategy, as we do here. However they strongly use the property that the free energy functional is bounded from below, which unfortunately does not hold true in our context, at least in the usual sense (see Lemma \ref{lem:Flowerbound}). This makes the analysis more difficult  and justifies {\em a priori} the weaker result which is obtained here.

We consider throughout this section $\e$ being positive but $\alpha$ being zero, since our aim here is to emphasize the parabolic character of the chemical equation. Then, our key strategy is to consider the parabolic-parabolic system as a perturbation of the parabolic-elliptic system. Doing that, the analysis of the time evolution of the second moment leads naturally to carefully estimate the time evolution of the corrected energy $\mathcal F[n,c]$ defined in \eqref{eq:corrected energy}. It will be shown that $\mathcal F[n,c]$  satisfies a master differential inequality (see \eqref{eq:master} below), generalization to the parabolic-parabolic Keller-Segel system of the differential inequality for $\mathcal F[n]$ obtained in the parabolic-elliptic case (see Remark \ref{rk:correctedenergy}). The master inequality together with the second moment growth rate estimate, give us the concentration result. Finally, the perturbation strategy could not work without the following two fundamental lemmas. 

\begin{lemma} [Free energy minimization]
Let $d\ge3$. Let $n$ be any nonnegative function in \break $(L^1\cap L^{2d/(d+2)})(\R^d)$, such that $\int_{\R^d}n(x)\log n(x)\, dx<\infty$. Let $\overline{c}:=E_d*n$. Then, the energy functional defined in \eqref{eq:energy} with $\alpha=0$, satisfies
\begin{equation}
\mathcal E[n,c]\ge \mathcal E[n,\overline c]=\int_{\mathbb{R}^d} n(x)\log n(x)\, dx-\f12 \int_{\R^d}n(x)\,\overline{c}(x)\, dx\ ,
\label{eq:energymin}
\end{equation}
for any $c$ such that $nc\in L^1(\R^d)$ and $|\nabla c|\in L^2(\R^d)$.
\label{lem:energymin}
\end{lemma}
%

\begin{proof}
In order to prove \eqref{eq:energymin} formally, it is sufficient to observe that
\begin{equation}
\mathcal E[n,c]-\mathcal E[n,\overline c]=\f12\int_{\mathbb{R}^d} |\nabla c(x)|^2\, dx-
\f12\int_{\mathbb{R}^d} |\nabla \overline c(x)|^2\, dx - \int_{\mathbb{R}^d}n(x)(c(x)-\overline c(x))\, dx\,.
\label{eq:energymin2}
\end{equation}
Since $-\Delta\overline c=n$, a simple integration by part in \eqref{eq:energymin2} gives us
\[
\mathcal E[n,c]-\mathcal E[n,\overline c]=\f12\int_{\mathbb{R}^d} |\nabla(c(x)-\overline c(x))|^2\, dx
\]
and the inequality in \eqref{eq:energymin} is proved, as well as the equality by applying  again an integration by part. Next, to rigorously justify the previous computation, we observe that $E_d\in L^{\f d{d-2}}_w(\R^d)$ and $|\nabla E_d|\in L^{\f d{d-1}}_w(\R^d)$. Therefore, we have on the one hand $\overline c\in L^{\f{2d}{d-2}}(\R^d)$, and on the other hand $|\nabla\overline c|\in L^{2}(\R^d)$  by the weak Young's inequality. Consequently, the integration by part in \eqref{eq:energymin2} can be performed (see for instance \cite{LiebLoss}).
\end{proof}
\begin{lemma} [Corrected energy lower bound]
Let $n$ be any nonnegative function in $(L^1\cap L^{\f d2})(\R^d)$ with finite second momentum. Let $c$ be such that $nc\in L^1(\R^d)$ and $|\nabla c|\in L^2(\R^d)$. With $M=\|n\|_{L^1(\R^d)}$ and the constant $B(d,M)$ defined in \eqref{eq:bdm}, the following lower bound for the corrected energy $\mathcal F[n,c]$ holds true
\begin{equation}
d(d-2)M\mathcal F[n,c]+B(d,M)+\f1{|\Sph^{d-1}|}M\,C_{HLS}(d,d-2)\|n\|_{L^{\f d2}(\R^d)}\ge 2dM\,.
\label{eq:Fbound}
\end{equation}
\label{lem:Flowerbound}
\end{lemma}
\begin{proof}
Let $\overline{c}=E_d*n$. Using the free energy minimization Lemma \ref{lem:energymin} and the expression for $B(d,M)$, we have
\begin{eqnarray}
d(d-2)M\mathcal F[n,c]+B(d,M)
&\ge& 2(d-2)\left[\f{dM}2\log I+\f{dM}2-M\log M-\f{dM}2\log\left(\f{dM}{2\pi}\right)\right]+2dM
\nonumber\\
& &\quad+\,2(d-2)\,\mathcal E[n,\overline c]
\nonumber\\ 
&=& 2(d-2)\,\left[\int_{\R^d}n\log n+\f{dM}2(\log I+1)-M\log M-\f{dM}2\log\left(\f{dM}{2\pi}\right)\right]
\nonumber\\
& &\quad+2dM-(d-2)\int_{\R^d}n\overline c\,.
\label{eq:Fbound2}
\end{eqnarray}
Since the bracket in the r.h.s. of \eqref{eq:Fbound2} is nonnegative by the entropy lower bound \eqref{est:entropy2}, inequality \eqref{eq:Fbound} follows from \eqref{eq:Fbound2} when applying to the quasi-stationary potential $\int_{\R^d}n\overline c$, the potential confinement inequality \eqref{est:potential}.

\end{proof}

Before proceeding in proving Theorem \ref{th:concentration}, let us observe that the constant $C(d)$ in \eqref{eq:concentration} is the same constant as in the master inequality \eqref{eq:master}. It derives from explicit computations in the proof of Theorem \ref{th:concentration}. Therefore, it can be quantified. Moreover, it is important to underline that the initial concentration condition \eqref{eq:concentration_condition} from one hand converges to the blow-up condition \eqref{2BU} as $\e\to0$ and from the other hand it implies that $\e$ is small with respect to $\|n_0\|_{L^{\f d2}(\R^d)}$. Indeed, applying the logarithmic function to both side of \eqref{eq:concentration_condition}, after rearranging the term, \eqref{eq:concentration_condition} reads as
\begin{equation}
d(d-2) M \mathcal F[n_0,c_0]+B(d,M)+d(d-2)M\,\e^\gamma<0\,.
\label{eq:concentration_condition_equiv}
\end{equation}
Then, the corrected energy lower bound   \eqref{eq:Fbound} gives 
\begin{equation}
d(d-2)\,\e^\gamma+2d<|\Sph^{d-1}|^{-1}\,C_{HLS}(d,d-2)\|n_0\|_{L^{\f d2}(\R^d)}\,.
\label{eq:small_epsilon}
\end{equation}
Obviously, inequality \eqref{eq:small_epsilon} does not imply that Theorem \ref{th:concentration} holds true only for small $\e>0$. Nevertheless, inequality \eqref{eq:concentration} gives a relevant concentration result only if $\e$ is small and $\|n_0\|_{L^{\f d2}(\R^d)}$ is upper bounded by the constant given in the r.h.s. of \eqref{eq:concentration}. One can try to obtain a concentration result for $n(t)$ when $\e$ is large, for exemple taking $\gamma>1$ or a more general function of $\e$ instead of $\e^\gamma$ into \eqref{eq:concentration_condition}. However, the only result obtained is that the $L^{d/2}$-norm of $n(t)$ stays large for a while when the $L^{d/2}$-norm of $n_0$ is large (see also Remark \ref{rk:concentration}). Finally, the maximal time of existence $T_{max}$ in Theorem \ref{th:concentration} is not requested to be finite and   we cannot exclude that the solution $(n,c)$ is global in time.

\begin{proof}[Proof of Theorem \ref{th:concentration}]
We shall proceed in several steps to conclude with a contradiction argument. The delicate and original points of the proof consist in the perturbation argument of the first step and in controlling the growth of the second moment $I(t)$, (no larger than $t\log t$), in the third step.\\

\noindent\emph{First step : the chemical quasi-stationary approximation.} We express the chemical $c$ as 
\begin{equation}
c=E_d*(n-\e\, \partial_t c)\,.
\label{eq:approxc}
\end{equation}
Indeed, applying the Fourier transform to the parabolic equation for $c$, we obtain 
$$
\widehat c\,(\xi)=\f1{4\pi^2|\xi|^2}\,(\,\widehat n-\e\, \partial_t\widehat c\,)(\xi)\,.
$$
Then, it suffices to apply the inverse Fourier transform and Theorem 5.9 or Corollary 5.10 in \cite{LiebLoss}, to obtain the identity \eqref{eq:approxc}. Consequently, the  gradient of $c$ can be written as follows:
\begin{equation}
\nabla c(x,t)=-\f{1}{|\mathbb{S}^{d-1}|}\int_{\R^d}\f{x-y}{|x-y|^d}(n(y,t)-\e\, \prt c(y,t))\, dy\,.
\label{eq:approxgradc}
\end{equation}
\vskip7pt

\noindent\emph{Second step : growth of $\mathcal F[n,c]$.}
Introducing the representation \eqref{eq:approxgradc} into the evolution equation \eqref{eq:evol I} for the second moment $I$ of $n$, we get after symmetrization
\begin{eqnarray}
\dfrac {d}{dt} I(t)&=&2dM-\frac{2}{|\mathbb{S}^{d-1}|}\iint_{\R^d\times\R^d}n(x,t)\,\f{x\cdot(x-y)}{|x-y|^d}\,n(y,t)\, dxdy
\nonumber\\
&&\quad\quad\quad
+\f{2\e}{|\mathbb{S}^{d-1}|}\iint_{\R^d\times\R^d}n(x,t)\,\f{x\cdot(x-y)}{|x-y|^d}\,\prt c(y,t)\, dxdy
\nonumber\\
\nonumber \\
& =& 2dM-(d-2)\int_{\R^d} n(t)\overline{c}(t)
\nonumber\\
&&\quad\quad\quad
+\f{2\e}{|\mathbb{S}^{d-1}|}\iint_{\R^d\times\R^d}n(x,t)\,\f{x\cdot(x-y)}{|x-y|^d}\,\prt c(y,t)\, dxdy\,.
\label{eq:I2para}
\end{eqnarray}
where $\overline c=E_d*n$. In order to control the second integral term in the r.h.s. of \eqref{eq:I2para}, we apply the general Hardy-Littlewod-Sobolev inequality \cite{LiebLoss}, and we use a suitable H\"{o}lder type inequality to get for every $\delta >0$,
\begin{eqnarray}
\!\!\!\!\!\iint_{\R^d\times\R^d}|x|n(x,t)\,\f1{|x-y|^{d-1}}\,|\prt c(y,t)|\, dxdy\!\!
&\le&\!\!\! C(d)\, \|x\, n(t,x)\|_{L^\f{2d}{d+2}(\R^d)}\, \|\partial_t c(t)\|_{L^2(\R^d)}
\nonumber\\
&\le&\!\!\! C(d)\,I^{\f12}(t)\|n(t)\|_{L^{\f d2}(\R^d)}^{\f 12}\|\prt c(t)\|_{L^2(\R^d)}
\nonumber\\
& \le&\!\!\! C(d)\!\!\left[\delta\,\|n(t)\|_{L^{\f d2}(\R^d)}+\delta^{-1}I(t){\|\prt c(t)\|_{L^2(\R^d)}^2}\right]\!\!,\  
\label{est:evolu1}
\end{eqnarray}
where the constant $C(d)$ above may change from line to line, the last being the one appearing in \eqref{eq:concentration}. Plugging estimation \eqref{est:evolu1} into \eqref{eq:I2para}, we obtain our first main estimate (the equivalent of \eqref{eq:I1} for $\varepsilon>0$)
\begin{equation}
\dfrac {d}{dt} I(t)\le 2dM-(d-2) \int_{\R^d} n(t)\overline{c}(t)+ \e\,\delta\,C(d)\,
\|n(t)\|_{L^{\f d2}(\R^d)}+\e\,\delta^{-1}C(d)\,I(t){\|\prt c(t)\|_{L^2(\R^d)}^2}\,.
\label{eq:I3para}
\end{equation}
Next, we use the free energy minimization Lemma \ref{lem:energymin} and the free energy dissipation equation \eqref{eq:energy dissipation}  into \eqref{eq:I3para} to have
\begin{equation}
\dfrac {d}{dt} I(t)\le 2dM+ 2(d-2)\left[\mathcal {E}[n,c](t) -\int_{\R^d} n(t)\log n(t)\right]
+\e\delta\,C(d)\,\|n(t)\|_{L^{\f d2}(\R^d)}-\delta^{-1}C(d)\,I(t)\,\dfrac {d}{dt} \mathcal{E}[n,c](t)\,. 
\nonumber
\end{equation}
Choosing $\delta:=\f{dM}2\,C(d)$, the entropy lower bound \eqref{est:entropy2} gives us
\begin{align}
\dfrac {d}{dt} I(t)\le\ &2dM+ 2(d-2)\,\mathcal {E}[n,c](t)+2(d-2)\left[\f{dM}2\log I+\f{dM}2-M\log M-\f{dM}2\log\left(\f{dM}{2\pi}\right)\right]
\nonumber \\
&+\f\e2\,dM\,C^2(d)\,\|n(t)\|_{L^{\f d2}(\R^d)}-\f2{dM}\,I(t)\,\dfrac {d}{dt} \mathcal{E}[n,c](t)\,,
\nonumber
\label{eq:I4para} 
\end{align}
i.e., after rearranging the terms, we get the following differential inequality for the corrected energy \eqref{eq:corrected energy}: 
\begin{equation}
I(t)\,\dfrac d{dt}\mathcal F[n,c](t)\le
d(d-2)M\,\mathcal F[n,c](t)+B(d,M)+\f\e2\,dM\,C^2(d)\,\|n(t)\|_{L^{\f d2}(\R^d)}\,.
 \label{eq:master}
 \end{equation}

\noindent\emph{Third step : moment's growth.}
Here we compute the evolution of $I(t)$ in an alternative way, involving the Fisher's information in the the energy dissipation equation \eqref{eq:energy dissipation}. As a consequence, we shall prove that the second  moment does not increase asymptotically faster than $2dM\,t\log t$ as long as $\|n(t)\|_{L^{\f d2}(\R^d)}$ stays bounded from above. Indeed,
\begin{eqnarray*}
\f{d}{dt} I(t) &=&  -2\int_{\R^d} n(x,t)\ x\cdot \nabla (\log n(x,t)-c(x,t))\, dx \\
&\leq& 2\left(\int_{\R^d} |x|^2n(x,t)\, dx\right)^{1/2}\left(\int_{\R^d} n(x,t)|\nabla (\log n(x,t)-c(x,t))|^{2}\,  dx\right)^{1/2}\, .
\end{eqnarray*}
Therefore,
\[ 
\f{d}{dt}I^{1/2}(t)\leq \left(\int_{\R^d} n(x,t)|\nabla (\log n(x,t)-c(x,t))|^{2}\,  dx\right)^{1/2}\, .
\]
Integrating the above inequality over $(0,t)$ we get
\[ 
I^{1/2}(t)\leq I^{1/2}(0)+\int_0^t  \left(\int_{\R^d} n(x,s)|\nabla (\log n(x,s)-c(x,s))|^{2}\,  dx\right)^{1/2}\, ds\, ,
\]
and 
\[I(t)\leq 2I(0)+2\,t\,\int_0^t\int_{\mathbb{R}^d} n(x,s)|\nabla (\log n(x,s)-c(x,s))|^2\, dxds\, .
\]
Using the  dissipation of the energy, we derive the  following pointwise estimate for $I(t)$ with respect to the corrected energy $\mathcal F[n,c](t)$
\begin{equation*}
I(t)\leq 2I(0)+2\,t\left(\mathcal E[n_0,c_0]- \mathcal E[n,c](t) \right )=
2I(0)+ 2\,t\,\mathcal E[n_0,c_0]-dMt\,\mathcal F[n,c](t)+dMt\,\log I(t)\, .
\end{equation*}
Finally, being $dMt\,\log I\le\f 12\,I+dMt\,\log(2dM\,t)$, we obtain 
$$
I(t)\leq 4I(0)+ 4\,t\,\mathcal E[n_0,c_0]-2dMt\,\mathcal F[n,c](t)+2dMt\,\log(2dM\,t)\,,
$$
i.e. the claimed behaviour for $I(t)$, thanks to the corrected energy lower bound \eqref{eq:Fbound},
\begin{equation}
I(t)\leq 4I(0)+ 4t\,\mathcal E[n_0,c_0]+2t\left[(d-2)^{-1}B(d,M)+\mu_d\,M\,C_{HLS}(d,d-2)\|n(t)\|_{L^{\f d2}(\R^d)}\right]
+2dMt\,\log(2dM\,t)\,.
\label{eq:momentgrowth}
\end{equation}

\noindent\emph{Fourth step : conclusion.} We now conclude showing by a contradiction argument that the concentration result \eqref{eq:concentration} holds true. Indeed,  comparing the master equation \eqref{eq:master} and the concentration condition \eqref{eq:concentration_condition} written as \eqref{eq:concentration_condition_equiv}, we claim that we can not have uniformly in time 
\begin{equation}
\f\e2\,dM\,C^2(d)\,\|n(t)\|_{L^{\f d2}(\R^d)}<d(d-2)M\,\e^\gamma\ .
\label{eq:contradiction}
\end{equation}
In order to proceed, we have to distinguish between two cases according to the fact that the initial density $n_0$ satisfies or not \eqref{eq:contradiction}. If $n_0$ does not satisfy \eqref{eq:contradiction}, the concentration result \eqref{eq:concentration} is obvious. If not we deduce  from the master equation \eqref{eq:master} and from \eqref{eq:concentration_condition_equiv} that the corrected energy is initially decreasing. Next, if \eqref{eq:contradiction} holds true uniformly in time, there exists $\delta>0$ such that $I(t)\frac d{dt}\mathcal F[n,c](t)<-\delta$, uniformly in time, i.e. $\mathcal F[n,c]$ remains decreasing for $t>0$. Moreover, plugging the upper bound for $\|n(t)\|_{L^{\f d2}(\R^d)}$ given by \eqref{eq:contradiction} into \eqref{eq:momentgrowth}, we get the estimate
\begin{equation}
\dfrac d{dt}\mathcal F[n,c](t)\leq -\dfrac{\delta}{4I(0)+C(\e,d,M,\mathcal E[n_0,c_0])\, t+2dM\,t\log t}\,.
\label{eq:non integrable in time}
\end{equation}
Since the right-hand-side of \eqref{eq:non integrable in time} is not integrable at infinity, but $\mathcal F[n,c](t)$ is bounded from below as soon as $\|n(t)\|_{L^{\f d2}(\R^d)}$ is bounded from above. We have obtained a contradiction which completes the proof of the theorem.
\end{proof}

\begin{remark} It is worth noticing that when  $n_0$ satisfies \eqref{eq:contradiction}, from \eqref{eq:small_epsilon} $\e$ results necessarily small, i.e.
\begin{equation}
\e<\f{2\,C_{HLS}(d,d-2)}{d|\Sph^{d-1}|\,C^2(d)}\,.
\label{eq:eps small}
\end{equation}
Moreover, the above upper bound for $\e$ is independent from the function of $\e$ choosen in the initial concentration condition \eqref{eq:concentration_condition}. On the other hand, if $n_0$ does not satisfy \eqref{eq:contradiction}, we have
$$
\|n_0\|_{L^{\f d2}(\R^d)}>\max\left\{d(d-2)|\Sph^{d-1}|C_{HLS}^{-1}(d,d-2)\,\e^\gamma\,;\,
\f{2(d-2)}{C^2(d)}\,\e^{\gamma-1}\right\}\,
$$
and the case ``$\e$ large'' is included here.
\label{rk:concentration}
\end{remark}
\begin{remark}[Concentration result of order $\e^{-1}$]
It is possible to obtain a concentration result of order $\e^{-1}$ under the ``critical'' hypothesis 
\begin{equation}
d(d-2)M\,\mathcal F[n_0,c_0]+B(d,M)+\f\e2\,dM\,C^2(d)\,\|n_0\|_{L^{\f d2}(\R^d)}<0\,,
\label{eq:critical concentration condition}
\end{equation}
that writes equivalently as 
$$
\int_{\R^d} |x|^2 n_0(x)\, dx< K_2(d)\
M^{1+ \f{2}d}\exp{\left(-\f2{dM}\mathcal E[n_0,c_0]\right)}
\exp\left(-\f\e{2(d-2)}C^2(d)\|n_0\|_{L^{\f d2}(\R^d)}\right)\,.
$$
Then, from the master equation \eqref{eq:master} we obtain that the corrected energy $\mathcal F[n,c]$ is initially decreasing. Proceeding as before by contradiction, we deduce the existence of  a time $T>0$ such that
$$
d(d-2)M\,\mathcal F[n_0,c_0]+B(d,M)+\f\e2\,dM\,C^2(d)\,\|n(T)\|_{L^{\f d2}(\R^d)}\geq 0\ .
$$
However, the critical concentration initial condition \eqref{eq:critical concentration condition} implies that $\e$ is necessarily small, i.e. $\e$ must satisfies \eqref{eq:eps small}.
\label{rk:concentration2}
\end{remark}
\begin{remark}
Theorem \ref{th:concentration} does not give any clue about the solution's behaviour after the first ``concentration time" (blow-up, persistence, dispersion ?). This is due to the fact that the master equation \eqref{eq:master} does not bring any information about the evolution of the second moment $I(t)$ and of $\|n(t)\|_{L^{\f d2}(\R^d)}$ after that time.
\end{remark}
\begin{remark}
It would be possible to extend the first blow-up criterion given in Proposition~\ref{prop:1BUcriterion} to the fully parabolic system as well. However the method is less natural and the final result is weaker. Thus we will not develop the arguments in this paper.
\label{rmk:concentration3}
\end{remark}
%

\section{A discrete model mimicking the parabolic-elliptic system}
\label{Sec:num}
In this section we aim to give a geometric intuition of what could be  the dynamics of the parabolic-elliptic Keller-Segel system in space dimension $d\geq 3$. For this purpose and for numerical convenience, first we shall replace the classical Keller-Segel system \eqref{KS} with $\eps=\alpha=0$, by a one-space dimensional variant having an analogous behaviour to \eqref{KS}. Energy and especially homogeneity considerations have guided our choice for the substitutive system as we shall see below. Next, we shall rephrase the one dimensional Keller-Segel system using the pseudo-inverse distribution function of the cell density $n$ (see \cite{BCC, GT2}). A numerical scheme following an idea of \cite{BCC} is discussed and a discrete dynamical system is then proposed, with the aim of providing a visualization tool for the dynamics of the Keller-Segel system.

Let us recall the action of the mass-preserving dilation $f(x)\to f^\lambda(x)=\lambda^{-d}\,f(\lambda^{-1}x)$ on the free energy \eqref{eq:energyPE} with $\alpha=0$, yet used in the subsection \ref{sec:complementarity},
\begin{equation}
\mathcal E [n^\lambda ]=  \boldsymbol{-dM\log\lambda} + \int_{\R^d} n(x,t) \log n(x,t) \, dx - \frac{\boldsymbol{\lambda^{2-d}}}2 \iint_{\R^d\times\R^d} n(x,t)E_d(x-y) n(y,t)\,dxdy\,.
\label{eq:lambdaenergy}
\end{equation}
We clearly see the different homogeneities of the two contributions composing the energy, namely the entropy and the potential. The first one has say an ``almost zero'' homogeneity (up to a logarithmic correction), while the second has the same degree of homogeneity as the kernel $E_d$.  These different homogeinities make the behaviour of the \eqref{KS} system more intricated in high dimension than in 2 dimension (see Remark \ref{rmk:2Dcase}). In order to reproduce the same high-dimensional behaviour in the simpler frame of the one space dimension, we consider from now on the following one parameter family of systems describing self-attracting diffusive particles
\begin{equation}
\left\{
\begin{array}{rcl}
\partial_t n& =& \partial_{xx}\, n - \partial_x(n\, \partial_x c)\, ,\\
c& =& K_\gamma *n\, ,
\end{array}
\right.
\label{KS-PE-1D}
\end{equation}
where the interaction kernel is given by $K_\gamma (x):=  \gamma^{-1}\,|x|^{-\gamma}$, with $\gamma\in(0,1)$. As for the \eqref{KS} system, the new system \eqref{KS-PE-1D} is equipped with the decreasing free  energy
\begin{equation} 
\mathcal{E}_\gamma[n](t) =  \int_{\R} n(x,t)\log n(x,t)\, dx - \dfrac{1}{2} \int_{\R} n(x,t)\, (K_\gamma* n)(x,t)\, dx\,,
\label{eq:energy 1D}
\end{equation}
satisfying
\[ 
\dfrac d{dt} \mathcal{E}_\gamma[n](t) = - \int_\R n(x,t) \left|\partial_x\left(\log n(x,t) - (K_\gamma * n)(x,t)\right)\right|^2\, dx\, . 
\]
The action of the one dimensional mass-preserving dilation $f(x)\to f^\lambda(x)=\lambda^{-1}\,f(\lambda^{-1}x)$ on the new energy \eqref{eq:energy 1D} is then given by
\begin{equation}
\mathcal E [n^\lambda](t)=  \boldsymbol{-M\log\lambda} + \int_{\R} n(x,t) \log n(x,t) \, dx - \frac{\boldsymbol{\lambda^{-\gamma}}}2 \int_{\R} n(x,t)(K_\gamma*n)(x,t)\,dx\,,
\label{eq:lambdaenergy1D}
\end{equation}
so reproducing \eqref{eq:lambdaenergy}.

\begin{remark}
The case $\gamma = -1$ is not included here due to integrability issues. It is too singular in dimension 1, although it would naturally correspond to the Poisson kernel in dimension $d = 3$ as far as homogeneity is concerned. The generalized model \eqref{KS-PE-1D} will be the subject of another paper \cite{CalvezCarrilloPrep} (not  restricted to dimension 1 only). 
\end{remark}
\begin{remark}[The two dimensional case]
Let us recall that in the two dimensional case the interaction kernel \eqref{eq:Ed} for the parabolic-elliptic \eqref{KS} system with $\alpha=0$ is of logarithmic type, i.e. $E_2(x)=-\f1{2\pi}\log|x|$, and the energy \eqref{eq:energyPE} reads as
\begin{equation}
\mathcal E[n](t)=\int_{\R^2} n(x,t)\log n(x,t)\,dx + \dfrac1{4\pi}\iint_{\R^2\times\R^2} n(x,t)\,\log|x-y|\,n(y,t)\,dx\,dy\, .
\label{eq:energyPE2D}
\end{equation}
Therefore, under the two dimensional mass-preserving dilation, the two energy contributions have the same  ``almost zero'' homogeneity (up to a logarithmic correction) and
\[
\mathcal E[n^\lambda]=2M\left(\f M{8\pi}-1\right)\log\lambda + \mathcal E[n]\,.
\]
The above property, in the gradient flow interpratation, directly implies blow-up of any density $n$ having super-critical mass $M>8\pi$ (see \cite{BCC,CalvezCarrilloPrep}). 
\label{rmk:2Dcase}
\end{remark}
\begin{remark}[The limiting case $\gamma=0$]
It is possible to mimic the above two dimensional case with the help of the one dimensional system \eqref{KS-PE-1D}, if one replace the kernel $K_\gamma$ with $\tilde{K}_\gamma (x)\!=\! \gamma^{-1}(|x|^{-\gamma}-1)$ and let $\gamma$ go to 0. Doing that, the limit interaction kernel is given by $\tilde K_0(x)= -\log|x|$ and the corresponding free energy satisfies
\begin{equation}
\mathcal{E}_0[n^\lambda ] =M\,\left( \dfrac{M}2 - 1 \right)\log\lambda+\mathcal{E}_0[n ]   \, ,
\label{eq:homogeneity}
\end{equation}
under the one dimensional mass-preserving dilation, exactly as for the classical two dimensional Keller-Segel system  \eqref{KS}. As a consequence, system \eqref{KS-PE-1D}, with this particular choice of kernel, yields a critical mass phenomenon, $M=2$ being the critical threshold, as we shall see below (see Remark \ref{rmk:gamma=0}).  This system has been considered previously in~\cite{BilerWoyczynski,CPS}.
\end{remark}
\begin{remark}[Corrected energy homogeneity] 
It is worth noticing that, still in the parabolic-elliptic case with $\alpha=0$, the corrected energy $\mathcal F[n]$ given in \eqref{eq:corrected energy}, has the following striking homogeneity structure: the action of the mass-preserving dilation on the term $\log\left(\int_{\R^d} |x|^2 n(x,t)\, dx\right)$ cancels with the logarithmic entropy's contribution, so that
\[
\mathcal F[n^\lambda]=\int_{\R^d} n(x,t)\log n(x,t) \, dx - \f{\boldsymbol{\lambda^{2-d}}}{dM}\iint_{\R^d\times\R^d} n(x,t)E_d(x-y) n(y,t)\,dxdy\,.
\]
\end{remark}
%
\subsection{Gradient flow interpretation of the one dimensional system}
\label{subsec:ref}
Let us define the pseudo-inverse distribution function of any nonnegative $L^1(\R)$ function $f$ having mass $M$, as the function $X(m)$ defined on the interval $(0,M)$ as follows
\[ 
X(m) := \inf\left\{ x\in\R \, \left|\, \int_{-\infty}^{x} f(y)\, dy \geq m\right.\right\}\,,\quad m\in (0,M)\, . 
\]
Then, for the pseudo-inverse distribution function $X=X(m,t)$ of the density $n$, the energy functional \eqref{eq:energy 1D} rewrites as 
\begin{equation}
\mathcal E_{\gamma}[n](t) = \mathcal G_\gamma[X](t) =  - \int_0^M\log \left(\partial_m X(m,t)\right)dm - \dfrac 1{2\gamma} \int_0^M\int_0^M |X(m,t)-X(m',t)|^{-\gamma}\, dm\, dm'\,,
\label{eq:energy JKO}
\end{equation}
while system \eqref{KS-PE-1D} rewrites as follows \cite{GT2}:
\begin{equation}
\label{eq:KS inverse} 
-\partial_t X(m,t) = \partial_m \left( \dfrac1{\partial_m X(m,t)} \right) + \int_0^M\mathrm{sign}\left(X(m,t) - X(m',t)\right)|X(m,t) - X(m',t) |^{-\gamma-1}\, dm' \,,
\end{equation}
Furthermore, the action of the one dimensional mass-preserving dilation on the pseudo-inverse distribution function simply consists in multiplying by $\lambda$, i.e. $X^\lambda(m) = \lambda X(m)$. Therefore,  the same different homogeneities of the two contributions in the new energy functional~\eqref{eq:energy JKO}  clearly appear again
\[ 
\mathcal G_\gamma[X^\lambda](t) = \boldsymbol{- M\log \lambda } - \int_0^M\log \left(\partial_m X(m,t)\right)\, dm  - \dfrac {\boldsymbol{ \lambda^{-\gamma} }}{2\gamma} \int_0^M\int_0^M |X(m,t)-X(m',t)|^{-\gamma}\, dm\, dm'\,.
\]
\begin{proposition}[Gradient flow interpretation]
The integro-differential equation \eqref{eq:KS inverse} is the gradient flow of the energy functional $\mathcal G_\gamma[X]$ for the Hilbertian structure on $L^2(0,M)$: 
\begin{equation}
\partial_t X = -\nabla_{L^2} \mathcal G_\gamma[X]\, . 
\label{eq:flowgradient}
\end{equation}
In addition, the $L^2(0,M)$ norm of the pseudo-inverse distribution function $X$ of $n$ is the second momentum of $n$:
\[ \|X(\cdot,t)\|_{L^2(0,M)} = \int_\R |x|^2n(x,t)\, dx\, . \]
\label{prop:gradflow}
\end{proposition}

This interpretation is due to Jordan, Kinderlehrer and Otto \cite{JKO,Otto}. We refer to the books \cite{Villani.OT,AGS} for a comprehensive presentation of this interpretation in any dimension of space. Here, for sake of completeness, we give a sketch of the proof of identity \eqref{eq:flowgradient}.

\begin{proof} To prove \eqref{eq:flowgradient}, we compute formally the variation of the functional $\mathcal G_\gamma[X]$ directly from \eqref{eq:energy JKO} to obtain
\begin{align*}
\mathcal G_\gamma [X + H] &= - \int_{0}^M \log\left( \partial_mX(m) + \partial_mH(m)\right)\, dm \\ & \qquad -  \dfrac 1{2\gamma} \int_0^M\int_0^M |X(m)-X(m') + H(m)-H(m')|^{-\gamma}\, dm\, dm'\\
& = \mathcal G_\gamma[X] -  \int_0^M  \dfrac{ \partial_mH(m)}{\partial_mX(m)}\, dm \\ & \qquad  +  \dfrac 1{2} \int_0^M\int_0^M \left|X(m)-X(m')\right|^{-\gamma}\dfrac{H(m)-H(m')}{X(m)-X(m')}\, dm\, dm' + \mathcal O(H)\\
& = \mathcal G_\gamma[X]  +  \int_0^M\partial_m\left( \dfrac{ 1}{\partial_mX(m)}\right)H(m)\, dm  \\ & \qquad + \int_0^M\int_0^M \left|X(m)-X(m')\right|^{-\gamma}\dfrac{H(m)}{X(m)-X(m')}\, dm\, dm' + \mathcal O(H)\,.
\end{align*}
\end{proof}
\begin{remark}[The critical mass phenomenon in the limiting case $\gamma=0$]
The energy functional $\mathcal G_0 [X]$, given by \eqref{eq:energy JKO} when replacing the kernel $K_\gamma$ with $\tilde K_0=-\log|x|$, satisfies
\[ 
\mathcal G_0 [\lambda X] = M\left(\dfrac{M}2 - 1\right)\log\lambda+\mathcal G_0 [X]\,,
\]
as one can deduce for exemple from \eqref{eq:homogeneity}. Differentiating the above relation with respect to $\lambda$ and evaluating the result at $\lambda = 1$ yields
\[ 
\left\langle X,\nabla \mathcal G_\gamma [X] \right\rangle_{L^2(0,M)} = M\left(\dfrac{M}2 - 1\right)\, . 
\]
On the other hand, $\langle X,\nabla \mathcal G_\gamma [X] \rangle_{L^2(0,M)}$ is nothing but the time derivative of $\|X(\cdot,t)\|_{L^2(0,M)}$ (up to a change of sign), as shown in Proposition \ref{prop:gradflow}. Therefore, we recover that the positive function $\|X(\cdot,t)\|_{L^2(0,M)}$ has a linear time decay when mass is super-critical, i.e. $M>2$.
\label{rmk:gamma=0}
\end{remark}
%


\subsection{The steepest descent scheme and the simplified dynamical system}
\label{subsec:steepest}
In \cite{JKO} the authors introduce a steepest descent numerical scheme for the linear Fokker-Planck equation. This coincides with the Euler time-implicit scheme when the equation is reformulated using the pseudo-inverse distribution function.  Those schemes possess the advantage of capturing important energy features (energy is decreasing and uniform  in time estimates can be derived). 

This strategy has been adapted to systems of  Keller-Segel type in \cite{BCC}, and particularly to the one-dimensional variant \eqref{KS-PE-1D} with the logarithmic interaction kernel $\tilde K_0$. It consists in discretizing the energy on a regular mesh in the mass space $(0,M)$, and then performing a gradient flow for the finite-dimensional functional obtained. As the numerical scheme is carefully built to preserve the gradient flow geometry of the original system, this procedure may give a clear intuition of the problem. In the sequel, we shall construct the same type of numerical scheme but starting from the energy functional \eqref{eq:energy JKO} with $\gamma>0$.

For sake of simplicity, we opt for a very rough discretization of the space $(0,M)$, namely we consider the following three points regular mesh
\[
m_1 = \dfrac M4 \, ,\quad m_2 = \dfrac M2\, , \quad m_3 = \dfrac {3}4M\, .
\]
  We next discretize the energy using finite differences. This writes, up to a constant factor $h$, as
\begin{equation}
{\bf G}_\gamma[X_1,X_2,X_3] = - \log(X_2 - X_1) - \log(X_3 - X_2)
- \dfrac h\gamma \left( (X_3 - X_2)^{-\gamma} + (X_3-X_1)^{-\gamma} + (X_2-X_1)^{-\gamma}\right)\,,
\label{eq:discrete energy}
\end{equation}
where $h = M/4$ is the space step, while the gradient flow of \eqref{eq:discrete energy} is given by the discrete dynamical system
\begin{equation}
\left\{
\begin{array}{l}
\dot X_1(t) = - \dfrac 1{X_2(t)-X_1(t)} +  h \left( (X_3(t) - X_1(t))^{-\gamma-1} + (X_2(t) - X_1(t))^{-\gamma-1} \right)\smallskip\\
\dot X_2(t) =  \dfrac 1{X_2(t)-X_1(t)} - \dfrac 1{X_3(t) - X_2(t)} + h \left( (X_3(t) - X_2(t))^{-\gamma-1} - (X_2(t) - X_1(t))^{-\gamma-1} \right)\smallskip\\
\dot X_3(t) =  \dfrac 1{X_3(t)-X_2(t)} - h \left( (X_3(t) - X_2(t))^{-\gamma-1} + (X_3(t) - X_1(t))^{-\gamma-1} \right)
\end{array}
\right.
\label{eq:discrete dynamical system}
\end{equation}

Let us observe that everything can be rewritten in terms of $u(t) = X_2(t)- X_1(t)$ and $v(t) = X_3(t)-X_2(t)$ because the center of mass $X_1(t)+X_2(t)+X_3(t)$ is conserved and assumed to be zero here, without loss of generality. Furthermore, we can reconstruct $X_1(t),X_2(t),X_3(t)$ from $u(t)$ and $v(t)$ by the following relations
\[ 
X_1(t) = - \dfrac13\left(2u(t) +v(t)\right)\, , \quad X_2(t) = \dfrac13\left(u(t) - v(t)\right)\, , \quad X_3(t) =  \dfrac13\left(u(t) + 2v(t)\right)\, . 
\]
The dynamics of the new variables $(u(t),v(t))$ are given by
\begin{equation} 
\left\{
\begin{array}{l}
\dot u(t) =  \dfrac 2{u(t)} - \dfrac 1{v(t)} + h \left( v(t)^{-\gamma-1} - 2 u(t)^{-\gamma-1} - (u(t)+v(t))^{-\gamma-1} \right)\smallskip\\
\dot v(t) = \dfrac 2{v(t)} - \dfrac 1{u(t)} + h \left( u(t)^{-\gamma-1} - 2 v(t)^{-\gamma-1} - (u(t)+v(t))^{-\gamma-1} \right)
\end{array}
\right.
\label{eq:BU toy model 1}
\end{equation}
 However, system \eqref{eq:BU toy model 1} is not a gradient flow in the updated variables $(u,v)$.
 
\begin{figure}
\begin{center}
\includegraphics[width = \linewidth,height = .8\textheight]{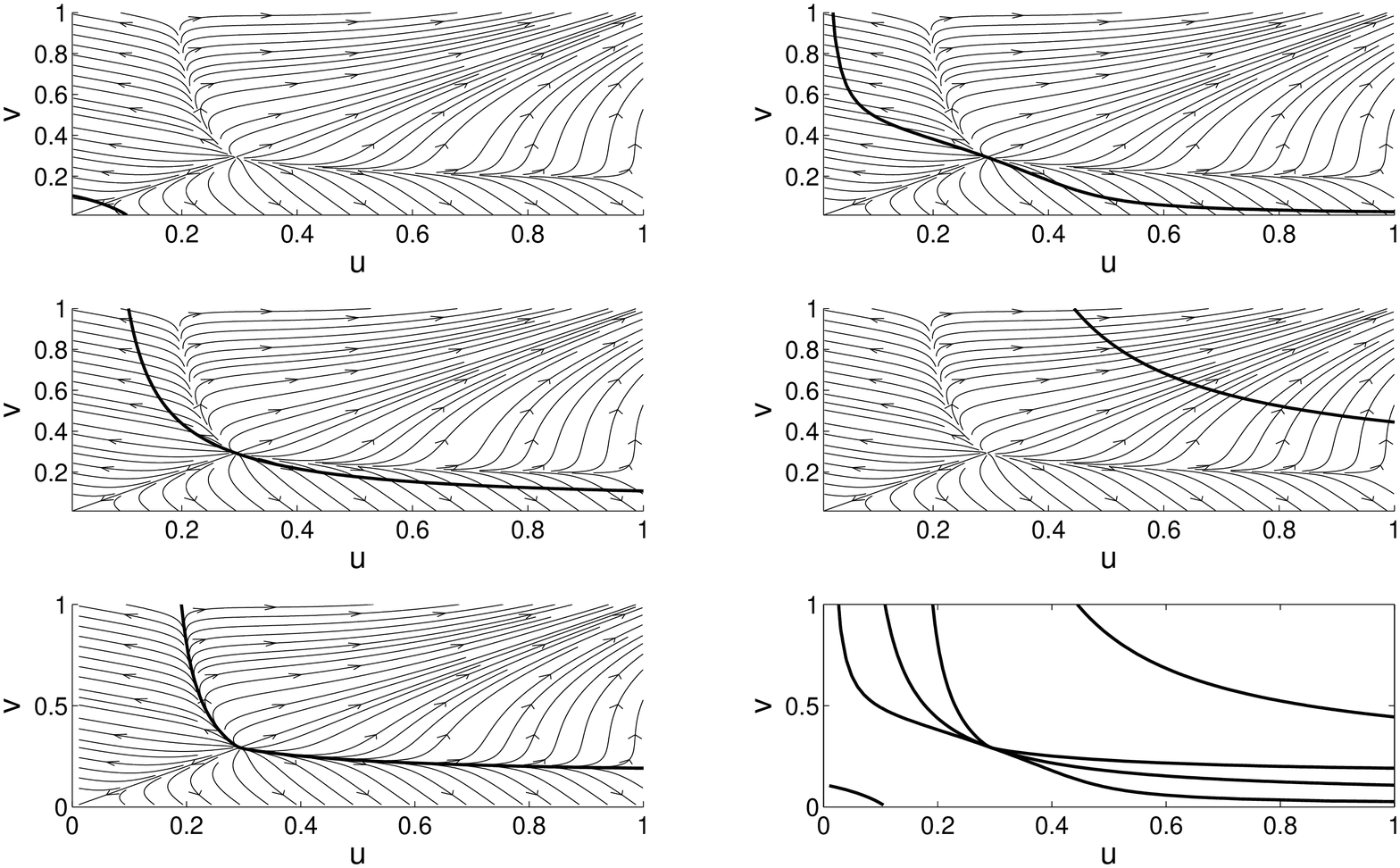}
\caption{\small Behaviour of the dynamical system \eqref{eq:BU toy model 1} in the phase plane $(u,v)$. ({\em Top}) The first  (\eqref{eq:criterion BU discrt 1}; {\em left}) resp. the second  (\eqref{eq:criterion BU discrt 2}; {\em right}) blow-up criterion  is figured as a bold line. ({\em Middle Left}) The line having equation ${\bf W}_\gamma[X] = 1$ (where the functional takes maximal values along radial rays) is plotted in bold. It is postulated to be an admissible criterion for blow-up but this has yet to be proved (see \eqref{eq:gauge definition} and the following discussion). (Middle Right) The global existence criterion is figured as a bold line. (Bottom Left) The unstable manifold (bold line) starting from the energy maximal point clearly separates the two basins of attraction. The derivation of its equation would provide a single criterion to distinguish between global existence and blow-up. (Bottom Right) Finally all the previous lines are plotted on the same figure for the sake of comparison.
}
\label{fig:BU non homogene}
\end{center}
\end{figure}

As an example we plot in Fig. \ref{fig:BU non homogene} the case $\gamma = 1/2$ and $M = 1.6$ (it would be similar for other choices of $M>0$ and $0<\gamma<1$). One clearly observes the separation between two opposite alternatives: ''blow-up'' when $u\to 0$ or $v\to 0$ and ''dispersion'' when both $u\to +\infty$ and $v\to +\infty$.

\subsection{Investigation of the landscape's geometry: blow-up vs. global existence}
\label{subsec:BUdiscrete}
In this section we shall derive two blow-up criteria and one global existence result for the dynamical system \eqref{eq:discrete dynamical system}. For that, it will be useful to rewrite the squared euclidean norm of $X = (X_1,X_2,X_3)$ in term of the updated variables $(u,v)$, 
\begin{equation}
|X|^2 = \dfrac23 \left(u^2 + v^2 + uv\right) = \dfrac13\left( u^2 + v^2 + (u+v)^2 \right)\,,
\end{equation}
as well as the functional ${\bf G}_\gamma[X]$
\begin{equation}
{\bf G}_\gamma[X]  = -\log u - \log v - \dfrac\chi\gamma \left( u^{-\gamma} + v^{-\gamma} + (u+v)^{-\gamma} \right)\, , \quad \chi = \dfrac{M}{4 }\, . 
\end{equation}

\begin{proposition}[Blow-up criteria]
\label{prop:criterion BU discrt 1} 
Assume that the initial point $X_0$ satisfies one of the two following criteria:
\begin{align}
& |X_0|^2< \left(\dfrac{3 \chi}2\right)^{2/\gamma} \, ,  \label{eq:criterion BU discrt 1}\\
& |X_0|^2\leq 2\exp\left( - {\bf G}_\gamma [X_0]-\dfrac2 \gamma\right)\, ,  \label{eq:criterion BU discrt 2}
\end{align}
then either $u(t)$ or $v(t)$ vanishes in finite time.
\end{proposition}
\begin{proof}
We follow carefully the evolution of the euclidean norm $|X(t)|$. The functional ${\bf G}_\gamma$ can be decomposed according to homogeneities:
\begin{equation}
{\bf G}_\gamma[X] = {\bf U}[X] - {\bf W}_\gamma[X] \, , \quad \left\{\begin{array}{l} {\bf U}[X] = - \log(X_2 - X_1) - \log(X_3 - X_2)\, ,\medskip \\ {\bf W}_\gamma[X] = \dfrac\chi\gamma \left( (X_3 - X_2)^{-\gamma} + (X_3-X_1)^{-\gamma} + (X_2-X_1)^{-\gamma}\right)\, .\end{array}\right.
\end{equation}
We derive the following Euler formula under mixed homogeneities assumptions:
\begin{equation}
\left\langle X,\nabla {\bf G}_\gamma [X] \right\rangle = \dfrac{d}{d\lambda}{\bf G}_\gamma[\lambda X]_{|_{\lambda = 1}}  = - 2  + \gamma{\bf W}_\gamma[X]\, .
\label{eq:connection with gauge function} 
\end{equation}
Therefore, the evolution of the euclidean norm under the gradient flow is driven by:
\begin{equation} 
\label{eq:toy L2 evolution}
\dfrac12\dfrac{d}{dt}|X(t)|^2 = - \left\langle X(t),\nabla {\bf G}_\gamma [X(t)] \right\rangle  =  2 - \gamma{\bf W}_\gamma[X(t)] \, .
\end{equation}
We use the following Jensen inequality based on the convexity of the function $(\cdot)^{-\gamma/2}$:
\begin{align}
\left(\dfrac13\left( u^2 + v^2 + (u+v)^2 \right)\right)^{-\gamma/2} & \leq \dfrac13\left( u^{-\gamma} + v^{-\gamma} + (u+v)^{-\gamma}\right)\, \label{eq:discrete ineq 1}  ,\\
|X|^{-\gamma} & \leq \dfrac\gamma{3\chi}{\bf W}_\gamma[X]  \, . \nonumber 
\end{align}
From \eqref{eq:toy L2 evolution} we obtain:
\begin{equation}
\dfrac12\dfrac{d}{dt}|X(t)|^2 \leq 2 - 3 \chi |X(t)|^{-\gamma}\, .
\end{equation}
Therefore, the norm necessarily vanishes in finite time if the initial data satisfies the following criterion: $|X_0|< (3 \chi/2)^{1/\gamma}$.

Starting from equation \eqref{eq:toy L2 evolution} we can plug the energy in the computation:
\begin{equation}
\dfrac12\dfrac{d}{dt}|X(t)|^2 \leq 2 + \gamma {\bf G}_\gamma [X_0] - \gamma {\bf U} [X(t)]\, .
\end{equation}
We use an alternative Jensen inequality, based on the concavity of $\log(\cdot)$:
\begin{align}
\log \left(\dfrac23\left(u^2+v^2+uv\right)\right)&\geq \log 2 + \dfrac13\left(2\log u + 2\log v + \log (uv)\right)\, , \label{eq:discrete ineq 2} \\
\log |X|^2 & \geq  \log 2 -  {\bf U}[X]\, . \nonumber
\end{align}
We obtain consequently:
\begin{equation}
\dfrac12\dfrac{d}{dt}|X(t)|^2 \leq 2 - \gamma \log 2+ \gamma {\bf G}_\gamma [X_0] +\gamma \log|X(t)|^2\, .
\end{equation}
Therefore, the norm necessarily vanishes in finite time if $|X_0|^2\leq 2\exp( - {\bf G}_\gamma [X_0]-2/\gamma)$.
\end{proof}

\begin{remark}[Comparison between the two criteria]
As opposed to Section \ref{Sec:Blow-up} the two criteria \eqref{eq:criterion BU discrt 1} and \eqref{eq:criterion BU discrt 2} are not complementary. Indeed the former implies the later (see Fig. \ref{fig:BU non homogene}). We shall prove this in the sequel. However there is a first point when looking for equality cases in resp. \eqref{eq:discrete ineq 1} and \eqref{eq:discrete ineq 2}. In fact, the former admits not equality case, whereas the later is an equality when $u = v$.

To prove that the first criterion \eqref{eq:criterion BU discrt 1} enhances the second \eqref{eq:criterion BU discrt 2}, we shall prove more generally that for all $X_0$,
\begin{equation} 
{\bf G}_\gamma [X_0] + \dfrac2\gamma  \leq \log 2 - \dfrac2\gamma \log\left(\dfrac{3   \chi}{2}\right) \, .  \label{eq:comparison discrete}
\end{equation}
The maximum of ${\bf G}_\gamma$ is achieved on the diagonal $\{u = v\}$. A simple computation yields for the extremal point: 
\[ u_*^{-\gamma} = \dfrac2{\chi  \left(2+2^{-\gamma}\right)}\, . \]
Plugging this into \eqref{eq:comparison discrete} we are reduced to prove:
\begin{equation}
\dfrac2\gamma \log\left(\dfrac3{ 2+2^{-\gamma} }\right)<\log2\, .\label{eq:comparison discrete 2}
\end{equation}
The function $\gamma\mapsto\log\left(2+2^{-\gamma}\right)$ being convex, we have for $\gamma\geq 0$:
\begin{equation*} \log\left(2+2^{-\gamma}\right) - \log 3  \geq - \dfrac\gamma3\log2   \geq - \dfrac\gamma2\log2 \, .
\end{equation*}
To conclude, the second criterion involving energy is clearly better in this case. This strongly uses the fact that the free energy is bounded from above (together with a precise evaluation of the maximal value), whereas this is no longer true in infinite-dimension. These computations are possible due to the simplicity of the toy model. 
\end{remark}

\begin{remark}[Derivation of a single criteria: open issue]
The separation between the basins of attraction of the two axes (blow-up) and of $\infty$ (global existence) is given by the unstable manifold of the energy critical point drawn on Fig. \ref{fig:BU non homogene}. We have failed in deriving an equation for this manifold.
\end{remark}

\begin{figure}
\begin{center}
\includegraphics[width = .44\linewidth]{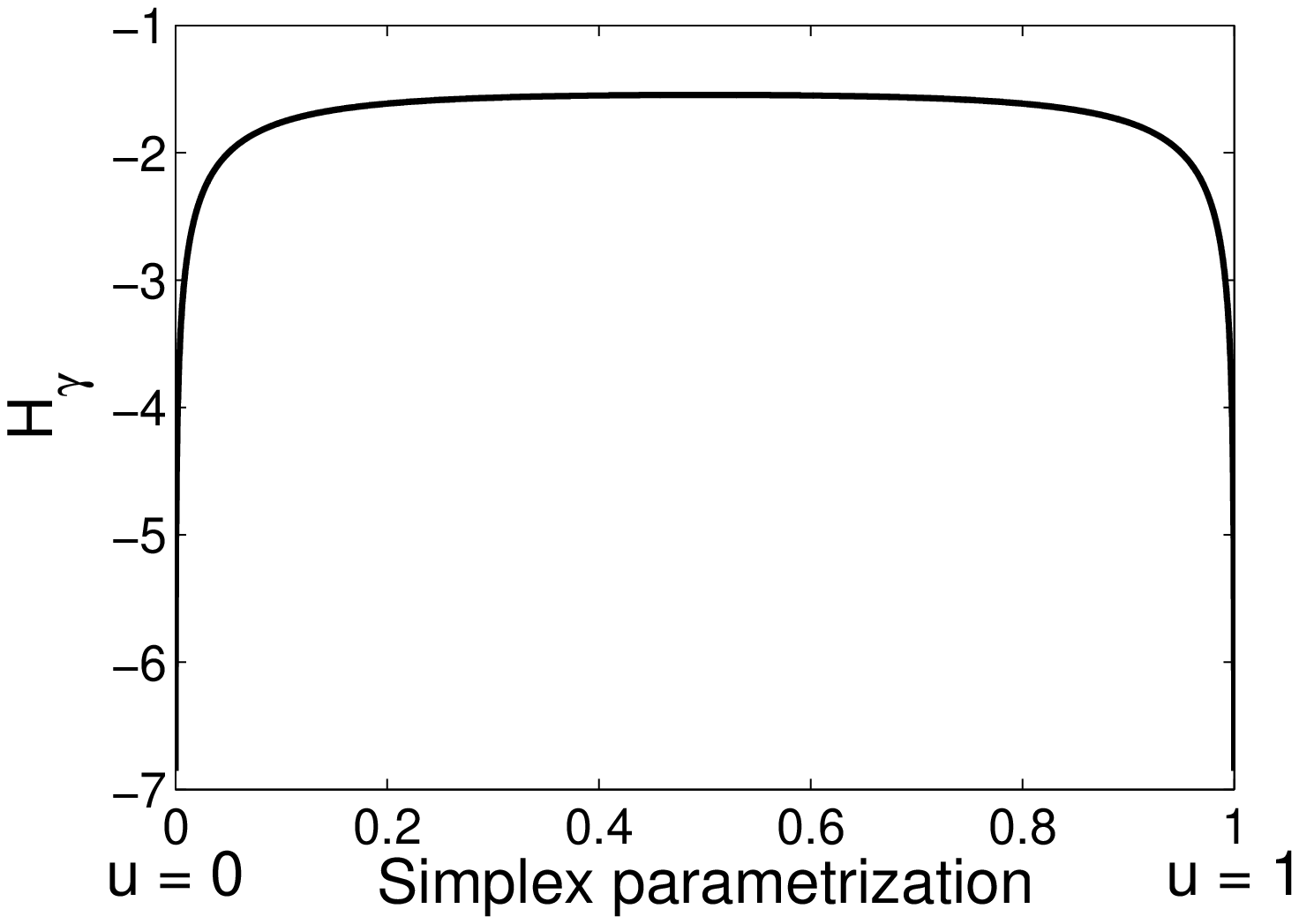}
\includegraphics[width = .52\linewidth]{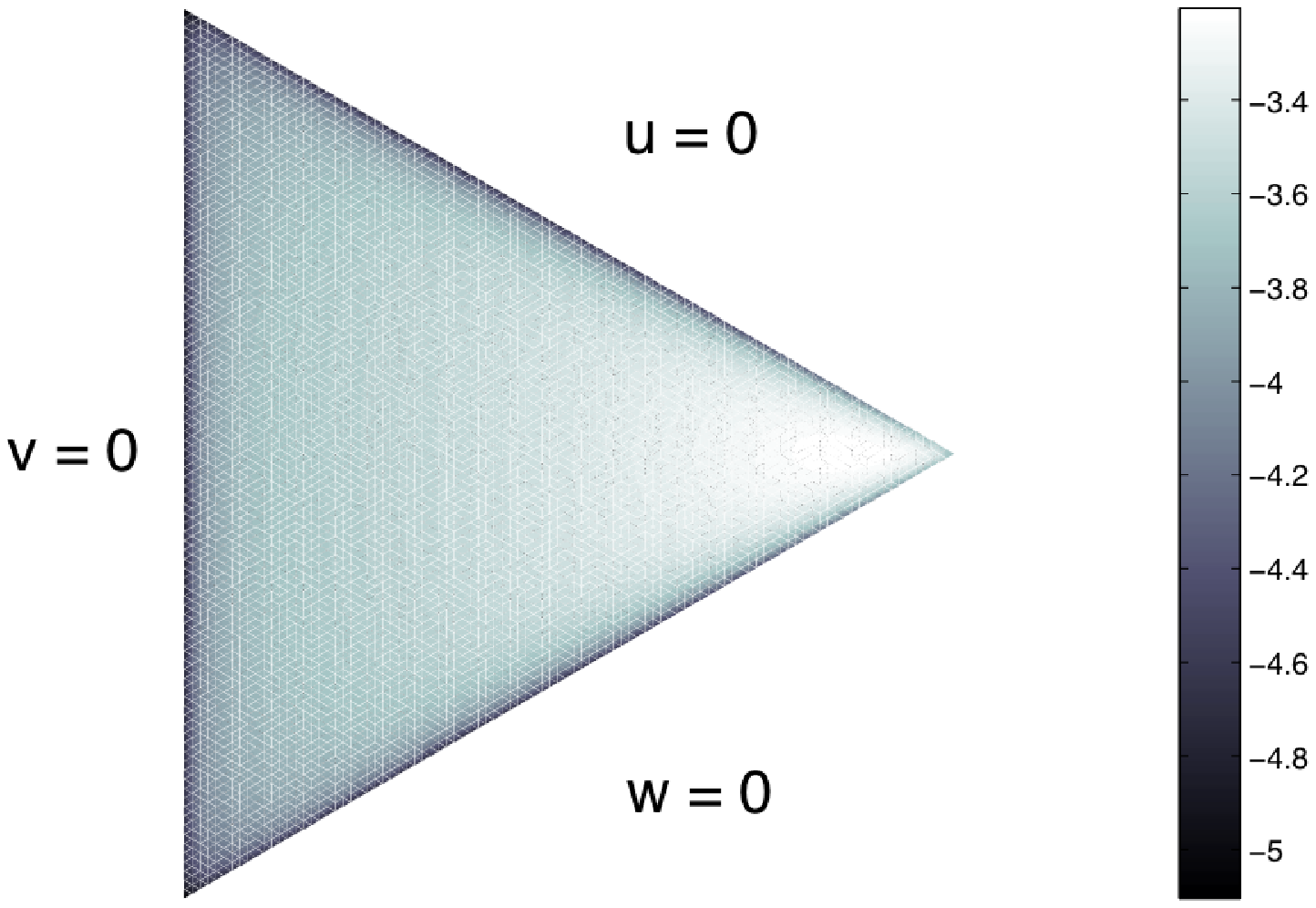}
\caption{\small Gauge function for (resp.) the three-points and four-points discretization of the energy functional. Recall that it is a function of homogeneity zero, which is defined by \eqref{eq:gauge definition}, where ${\bf U}[X] = -\log u - \log v $ (three points) or ${\bf U}[X] = -\log u - \log v - \log w $ (four points), and ${\bf W}_\gamma[X]  = \chi\gamma^{-1} \left( u^{-\gamma} + v^{-\gamma} + (u+v)^{-\gamma} \right)  $ (three points), and ${\bf W}_\gamma[X]  = \chi\gamma^{-1} \left( u^{-\gamma} + v^{-\gamma} + w^{-\gamma} + (u+v)^{-\gamma} + (v+w)^{-\gamma} + (u+v+w)^{-\gamma}\right)$ (four points). The maximal value is plotted on the simplex (resp. $u+v = 1$, $u+v+w = 1$) because it is a functional of homogeneity zero. 
}
\label{fig:gauge}
\end{center}
\end{figure}

\begin{remark}[Towards a better description of the landscape's geometry]
We introduce below another correction of the free energy. Notice that although it comes out very naturally, it is still unclear how to benefit from it appropriately. We aim to gain some insights concerning the geometry of the functional landscape, by means of homogeneity arguments. We decompose the free energy as previously:  ${\bf G}_\gamma[X] = {\bf U}[X] - {\bf W}_\gamma[X]$ where ${\bf U}$ is almost zero-homogeneous (up to a logarithmic correction) and ${\bf W}_\gamma[X]$ is $(-\gamma)$-homogeneous. We define the maximal value along rays:
\begin{equation}\label{eq:gauge definition}
{\bf H}_\gamma[X] = \max_{\lambda>0} {\bf G}_\gamma[\lambda X] = {\bf U}[X] - \dfrac2\gamma\left( \log\left(\dfrac{\gamma{\bf W}_\gamma[X]}2\right) + 1\right)\, .
\end{equation}
This function is zero-homogeneous by definition. Moreover, a simple computation shows that the set of extremal points: $\{ \lambda^* X\, |\, \lambda^* = {\rm argmax_\lambda}\,{\bf G}_\gamma[\lambda X]   \}$  admits the following equation: $ {\bf W}_\gamma[X] = 2/\gamma $ (cf. \eqref{eq:connection with gauge function}). Observe that the maximal value coincides with ${\bf U}[X] - 2/\gamma$ precisely on this curve (hypersurface in higher dimension). 

A reasonable claim would be that the gradient flow can not overpass this "maximal line", so that the later splits the phase space into two regions, and blow-up is guaranteed whenever we start initially on one side of this line. It turns out that this "maximal line" does not coincide with the unstable manifold which is the clear separation between global existence and blow-up (see Fig. \ref{fig:BU non homogene}). 

It is also worth noticing that the set $\{{\bf W}_\gamma[X]>2/\gamma\}$  coincides exactly with the area of the phase space where the second momentum is initially decreasing
\eqref{eq:toy L2 evolution}. The connection obviously derives from the definition of ${\bf H}_\gamma[X]$. 
\end{remark}

  To conclude we give a global existence result which is somehow analogous to Theorem \ref{the:localexistence}. This criterion is optimal in the sense that the asymptotes coincide with the  unstable manifold (not shown).
\begin{proposition}[Global existence]
If the initial data satisfies: 
\begin{equation}\chi   \left(u_0^{-\gamma} + v_0^{-\gamma}\right)<1\, ,
\label{eq:discrete GE crit}
\end{equation}
then the solution is global.
\end{proposition}
\begin{proof} 
To prove global existence in the line of Section \ref{Sec:LGexistence} we shall estimate some suitable quantity preventing $u(t)$ and $v(t)$ from vanishing. Namely, we investigate the evolution of $L(t) = \gamma^{-1}\left(u(t)^{-\gamma} + v(t)^{-\gamma}\right) $. We have:
\begin{align*}
\dfrac d{dt} L(t) & = - \left[ 2 u(t)^{-\gamma-2} + 2 v(t)^{-\gamma-2} - v(t)^{-1} u(t)^{-\gamma - 1} - u(t)^{-1} v(t)^{-\gamma - 1} \right] \\ & \qquad+ \chi  \left[ 2  u(t)^{-2\gamma-2} + 2  v(t)^{-2\gamma-2}  - 2 v(t)^{-\gamma-1} u(t)^{-\gamma-1} + (u(t)+v(t))^{-\gamma-1} \left(u(t)^{-\gamma-1} + v(t)^{-\gamma-1}\right) \right]  \, .
\end{align*}
We claim that the following inequality \eqref{eq:discrete GNS} holds true. This enables to conclude that the quantity $L(t)$ is nonincreasing in time if \eqref{eq:discrete GE crit} is verified. 
\begin{multline}
 2  u ^{-2\gamma-2} + 2  v ^{-2\gamma-2}  - 2 v ^{-\gamma-1} u ^{-\gamma-1} + (u +v )^{-\gamma-1} \left(u ^{-\gamma-1} + v ^{-\gamma-1}\right) \\\leq \left[ u^{-\gamma} + v^{-\gamma}\right]  \left[ 2 u ^{-\gamma-2} + 2 v ^{-\gamma-2} - v ^{-1} u ^{-\gamma - 1} - u ^{-1} v ^{-\gamma - 1}   \right]\, . 
\label{eq:discrete GNS}  \end{multline}

This is obvious when homogeneities coincides ($\gamma = 0$).
Due to homogeneity, \eqref{eq:discrete GNS} reduces to the following inequality for $U = u/v$:
\begin{equation}
 2  U ^{-2\gamma-2} + 2    - 2  U ^{-\gamma-1} + (U + 1 )^{-\gamma-1} \left(U ^{-\gamma-1} +1\right) \\\leq \left[ U^{-\gamma} + 1\right]  \left[ 2 U^{-\gamma-2} + 2   -  U ^{-\gamma - 1} - U ^{-1} \right]\, . 
\label{eq:reduced discrete GNS}
\end{equation}
We restrict to the case $U\leq1$ without loss of generality. We bound the delicate term $(U+1)^{-\gamma-1}$ by $1$. We end up with:
\[ 2 U^{-2\gamma-2} + 3 -  U^{-\gamma-1} \, , \]
on the one hand, and 
\[    2U^{-2\gamma-2} + 2 + 2 U^{-\gamma-2} + 2U^{-\gamma} - U^{-2\gamma-1} - U^{-1} - 2 U^{-\gamma-1} \, , \]
on the other hand.
Thus, we are reduced to prove that:
\begin{align} 
1 + U^{-2\gamma-1} + U^{-1} + U^{-\gamma-1} &\leq 2 U^{-\gamma-2} + 2 U^{-\gamma}\, , \\
1 + U^{-\gamma - 1} & \leq  U^{-\gamma} + U^{-\gamma - 2} +  U^{-1}\left( U^{-\gamma - 1} - 1\right)\left( 1 - U^{1 - \gamma} \right)  
\end{align}
the later being clear term by term because $\gamma\in (0,1)$ and $U\leq 1$.

\end{proof}

\section{Concluding remarks and open problems}
\label{Sec:OP}
In this paper we have investigated existence and blow-up issues concerning the Keller-Segel system \eqref{KS} in high space dimension (namely $d\geq 3$). Concerning the blow-up, we have analyzed two sufficient criteria for the parabolic-elliptic Keller-Segel system, and we have discussed complementarity between them. We have proposed a possible extension of one of  these criteria to the parabolic-parabolic case. Although we cannot still exhibit blow-up in that case, we are able to prove that the ``critical'' $L^{d/2}-$norm of the density $n$ reaches high level after some finite time. We have also proposed a discrete model for visualizing the dynamics of the parabolic-elliptic Keller-Segel model in finite dimension. All this analysis highlights the complexity of the high dimensional Keller-Segel system. Although in dimension 2 the existence and blow-up issues, in a simplified scheme, reduce to the mass being smaller or greater than $8\pi$, here it appears that the separation between global existence and blow-up is not so easy to describe, even for the discrete model.

We aim to conclude the paper listing below some problems which are still open up to our best knowledge.
\begin{enumerate}[(a)]
\item There is a gap between the global existence criterion \eqref{n_0condSob} and any blow-up criteria \eqref{1BU} and~\eqref{2BU}. Indeed, both of them imply the inequality
\begin{equation}
\iint_{\R^d\times\R^d}n_0(x)\f{1}{|x-y|^{d-2}}n_0(y)\,dxdy>2dM|\Sph^{d-1}|\,.
\label{eq:generalBU}
\end{equation}
which implies 
$$
\|n_0\|_{L^{\f d2}(\R^d)}>\f{2d|\Sph^{d-1}|}{C_{HLS}(d,d-2)}
=\f{2d}{(d-2)\,C^2_S(d)}>\f8{d\,C^2_S(d)}\,.
$$
Let us observe that \eqref{eq:generalBU} reduces to $M>8\pi$ in dimension 2.
\item The two blow-up criteria share the following feature: they imply that $\f d{dt}I(t)_{|_{t=0}}<0$, that is nothing else than \eqref{eq:generalBU}. Would it be possible to prove that the initial decreasing behaviour of the second moment is sufficient for blow-up ? This seems to be the case on the finite dimensional model. However this does not give an exhaustive description of the set of initial data for which blow-up occurs.
\item It has been shown that none of the two blow-up criteria implies the other. However in the discrete model, the criterion involving energy is clearly better. Is it possible to show that the energy criterion \eqref{2BU}, which is more convenient to deal with, is also better in some sense for the \eqref{KS} system ?
\item Are we able to derive the equation for the separation line (unstable manifold) in Figure \ref{fig:BU non homogene}? Observe that such an unstable manifold suggests the existence of a stationary state, at least in the discrete model for which the free energy is bounded from above. It is a challenging issue to understand better the geometry of the energy landscape for the parabolic-elliptic Keller-Segel system.
\end{enumerate}
\bigskip
{\bf Acknowledgement}. This paper has been partially prepared during the visit of the three authors to the Centre de Recerca Matem\`atica, Universitat Aut\`onoma de Barcelona, for the special semester ``Mathematical Biology: Modelling and Differential Equations''. They would all like to express their gratitude for the invitation and the kind hospitality. MAE has been supported by the Marie Curie Actions of the European Commission in the frame of the DEASE project (MEST-CT-2005-021122).

%
%


\noindent ${^a}$ Unit\'e de Math\'ematiques Pures et Appliqu\'ees,  CNRS UMR 5669 \\
\'Ecole Normale Sup\'erieure de Lyon, \\
46 all\'ee d'Italie, F~69364 Lyon cedex 07, France \\
e-mail: vincent.calvez@umpa.ens-lyon.fr
\bigskip

\noindent ${^b}$ D\'epartement de Math\'ematiques, \\
Universit\'e d'Evry Val d'Essonne, \\
Rue du P\`ere Jarlan, F~91025 Evry Cedex, France \\
e-mail: lucilla.corrias@univ-evry.fr
\bigskip

\noindent ${^c}$ 
Wolfgang Pauli Institute, \\ 
c/o Fak. f. Mathematik,
Nordbergstr. 15 A-1090 Wien, Austria\\
e-mail: mohamed.abderrahman.ebde@univie.ac.at


\begin{thebibliography}{99}
%
\bibitem{AGS} L. Ambrosio, N. Gigli and G. Savar\'e, Gradient flows in metric spaces and in the space of probability measures, Lectures in Mathematics, ETH ZŸrich. Birkh\"auser Verlag, Basel, 2005.
%
\bibitem{AU}
J.-P. Aubin, Un th\'eor\`eme de compacit\'e, {\em C. R. Acad. Sci. Paris} {\bf 256} (1963) 5042--5044.
\bibitem{Beckner} W. Beckner, Sharp Sobolev inequalities on the sphere and the Moser-Trudinger inequality,  {\em Ann. of Math.} {\bf 138}  (1993)  213--242.
\bibitem {B95} P. Biler,
Existence and nonexistence of solutions for a model of gravitational interaction of particle III,
{\em Colloq. Math.}  {\bf 68} (1995) 229--239.
\bibitem{BCD} P. Biler, L. Corrias and J. Dolbeault, Large mass self-similar solutions of the parabolic-parabolic Keller--Segel model of chemotaxis, Preprint arXiv 2009, \texttt {http://arxiv.org/abs/0908.4493v1}.
\bibitem{BilerWoyczynski}
P. Biler and W.A. Woyczy\'nski, Global and exploding solutions for 
nonlocal quadratic evolution problems, {\em SIAM J. Appl. Math.} {\bf 59} (1998)
845--869. 
\bibitem{BCC} A. Blanchet, V. Calvez and J.A. Carrillo, Convergence of the mass-transport steepest descent scheme for the subcritical Patlak-Keller-Segel model.  {\em SIAM J. Numer. Anal.}  {\bf 46}  (2008) 691--721.
\bibitem{BCL} A. Blanchet, J.A. Carrillo and Ph. Lauren\c{c}ot,  Critical mass for a Patlak-Keller-Segel model with degenerate diffusion in higher dimensions, {\em Calc. Var. Partial Differential Equations} {\bf  35}  (2009) 133--168.
\bibitem{BCM} A. Blanchet, J.A. Carrillo and N. Masmoudi, 
Infinite time aggregation for the critical Patlak-Keller-Segel model in $\R^2$, 
\newblock {\em Comm. Pure Appl. Math.} {\bf 61} (2008) 1449--1480.
\bibitem{BDP} A. Blanchet, J. Dolbeault and B. Perthame,
Two-dimensional Keller-Segel model: optimal critical mass and qualitative properties of the solutions,
\newblock {\em Electron. J. Diff. Eqns.}  {\bf 44} (2006)  1--33.
\bibitem{BournaveasC} N. Bournaveas and V. Calvez, Critical mass phenomenon for a chemotaxis kinetic model with spherically symmetric initial data,
{\em Annales de l'Institut Henri Poincare (C) Non Linear Analysis} {\bf 26} (2009) 1871--1895.
\bibitem{CalvezCarrilloPrep} V. Calvez and J.A. Carrillo, In preparation.
\bibitem{CC08} V. Calvez and L. Corrias,
The parabolic-parabolic Keller-Segel model in $\R^2$,
\newblock {\em Commun. Math. Sci.} {\bf 6} (2008) 417--447.
\bibitem{CPS}
V. Calvez, B. Perthame and M. Sharifi~tabar,
Modified {K}eller-{S}egel system and critical mass for the log interaction
  kernel,
in Nonlinear partial differential equations and related analysis, vol. 429 of {\em Contemp. Math.}, Amer. Math. Soc., Providence, RI, 2007.
\bibitem{CarlenLoss}
E.~Carlen and M.~Loss,
Competing symmetries, the logarithmic {HLS} inequality and {O}nofri's inequality on ${\Sph}\sp n$,
\newblock {\em Geom. Funct. Anal.} {\bf 2} (1992)  90--104.
\bibitem {CP81} S. Childress and J. K. Percus,
Nonlinear aspects of chemotaxis,
\newblock  {\em Math. Biosci.}  {\bf 56} (1981)  217--237.
\bibitem{CL} T. Cie\'slak and P. Lauren\c{c}ot,
Finite time blow-up for a one-dimensional quasilinear parabolic-parabolic chemotaxis system, Preprint arXiv 2009, 
\texttt{http://arxiv.org/abs/0810.3369}.
\bibitem {CP_CRAS06} L.  Corrias and B. Perthame,
Critical space for the parabolic-parabolic Keller-Segel model in $\R^d$,
\newblock {\em C.  R.  Acad.  Sci.  Paris}, Ser. I  {\bf 342} (2006) 745--750.
\bibitem {CP08} L.  Corrias and B. Perthame,
Asymptotic decay for the solutions of the parabolic-parabolic Keller-Segel chemotaxis system in critical spaces, 
\newblock{\em Math. Comp.  Model.} {\bf 47} (2008) 755--764.
\bibitem {CPZ04} L. Corrias, B. Perthame and H.  Zaag,
Global solutions of some chemotaxis and angiogenesis systems in high space dimensions,
\newblock  {\em Milano J. of Math.}  {\bf 72} (2004)  1--29.
\bibitem{DelPinoDolbeault} 
M. Del Pino and J. Dolbeault, Best constants for Gagliardo-Nirenberg inequalities and applications to nonlinear diffusions, {\em J. Math. Pures Appl.} {\bf 81} (2002) 847--875.
\bibitem{DolbeaultSchmeiser}
J. Dolbeault, C. Schmeiser, The two-dimensional Keller-Segel model after blow-up, {\em Discrete and Continuous Dynamical Systems} {\bf 25} (2009) 109--121. 
\bibitem{Gajewski98}
H.~Gajewski and K.~Zacharias,
Global behavior of a reaction-diffusion system modelling chemotaxis,
{\em Math. Nachr.} {\bf 195} (1998) 77--114.
\bibitem{G} R. T. Glassey,
On the blowing up of solutions to the Cauchy problem for nonlinear Schr\"{o}dinger equations,
\newblock{\em J. Math.Phys.} {\bf 18} (1977) 1794--1797.
\bibitem{GlasseyBook} 
R.T. Glassey, {\em The Cauchy problem in kinetic theory}. Society for Industrial and Applied Mathematics (SIAM), Philadelphia, PA, 1996.
\bibitem{GT2}
 L. Gosse and G. Toscani, Lagrangian numerical
approximations to one-dimensional convolution-diffusion
equations, {\em SIAM J. Sci. Comput.} {\bf 28} (2006) 1203--1227.
\bibitem {HP09} T. Hillen and K. Painter,
A user's guide to PDE models for chemotaxis,
\newblock  {\em J. Math. Biol.}  {\bf 58} (2009)  183--217.
\bibitem{Horstmann.Wang01}
D. Horstmann and G. Wang, Blow-up in a chemotaxis model without symmetry assumptions, {\em European J. Appl. Math.} {\bf 12} (2001) 159--177. 
\bibitem{Horstmann02} D. Horstmann, 
On the existence of radially symmetric blow-up solutions for the Keller-Segel model, {\em J. Math. Biol.} {\bf 44}  (2002) 463--478.
\bibitem{Horstmann03}
D. Horstmann, From 1970 until present : {the Keller-Segel}
model in chemotaxis and its consequences. {I}., {\em Jahresber.
Deutsch. Math.-Verein.} {\bf 105} (2003) 103--165.
%
\bibitem{HW05}
D. Horstmann and M. Winkler, Boundedness vs. blow-up in a chemotaxis system, {\em J. Differential Equations} {\bf 215}  (2005) 52--107.
\bibitem {JL} W. J\"ager and S. Luckhaus,
On explosions of solutions to a system of partial differential equations modeling  chemotaxis,
{\em Trans. Amer. Math. Soc.} {\bf  239} (1992) 819--821.
\bibitem{JKO} R. Jordan, D. Kinderlehrer and F. Otto, The variational formulation of the Fokker-Planck equation, {\em SIAM J. Math. Anal.} {\bf 29}  (1998) 1--17.
\bibitem {KS70} E.F. Keller and L.A. Segel,
Initiation of slime mold aggregation viewed as an instability,
{\em J. Theor. Biol.} {\bf 26} (1970) 399--415.
\bibitem{KS71} E.F. Keller and L.A. Segel, 
Model for chemotaxis, {\em  J. Theor. Biol.} {\bf  30} (1971) 225--234. 
\bibitem{KozoSugy} H. Kozono and Y. Sugiyama, 
Strong solutions to the Keller-Segel system with the weak $L^{\f n2}$ initial data and its application to the blow-up rate, {\em  Preprint}. 
\bibitem{KozoSugy08} H. Kozono and Y. Sugiyama, 
The Keller-Segel system of parabolic-parabolic type with initial data in weak $L^{\f n2}(\R^n)$ and its applications to Self-similar solutions, {\em Ind. Univ. math. J.} {\bf 57} (2008) 1467--1499. 
\bibitem{Lieb}
E.H. Lieb, Sharp constants in the Hardy-Littlewood-Sobolev and related inequalities, {\em Ann. of Math.} {\bf 118}  (1983) 349--374.
%
\bibitem{LiebLoss} E.H. Lieb and M. Loss,
\newblock {\em  Analysis}, second edition, Graduate Studies in Mathematics, vol. 14, American Mathematical Society (2001).
\bibitem{Nagai95}T.  Nagai, Blow-up of radially symmetric solutions to a chemotaxis system, {\em   Adv. Math. Sci. Appl.}  5  (1995) 581--601.
\bibitem{NSY} T. Nagai, T. Senba and K. Yoshida,
Application of the Trudinger-Moser inequality to a parabolic system of chemotaxis,
{\em  Funk. Ekv.} {\bf  40} (1997) 411--433.
\bibitem{NaitoSY} Y. Naito, T. Suzuki, K. Yoshida, 
Self-similar solutions to a parabolic system modeling chemotaxis. 
{\em J. Differential Equations} {\bf184} (2002) 386--421
\bibitem{Nanjundiah}
V. Nanjundiah, Chemotaxis, signal relaying and
aggregation morphology, {\em J. Theor. Biol.} {\bf 42} (1973) 63--105.
\bibitem{Otto} F. Otto, The geometry of dissipative evolution equations: the porous medium equation, {\em  Comm. Partial Differential Equations} {\bf 26}  (2001) 101--174.
\bibitem{P} B. Perthame,
{\em Transport Equation in Biology},
Frontiers in Mathematics, Birkh\"auser, 2007
\bibitem{VelazquezA}
J.J.L. Vel\'azquez, Point dynamics in a singular limit of the Keller-Segel model. {I}. Motion of the concentration regions, {\em SIAM J. Appl. Math.} {\bf 64} (2004) 1198--1223.
\bibitem{VelazquezB}
J.J.L. Vel\'azquez, Point dynamics in a singular limit of the Keller-Segel model. {II}. Formation of the concentration regions, {\em SIAM J. Appl. Math.} {\bf 64}  (2004) 1224--1248.
\bibitem{Villani.OT} C. Villani, {\em Topics in optimal transportation}. Graduate Studies in Mathematics, 58. American Mathematical Society, Providence, RI, 2003.
\bibitem{W} M. I. Weinstein,
Nonlinear Schr\"{odinger equations and sharp interpolation estimates},
\newblock {\em Commun. Math. Phys.} {\bf 87} (1983) 567--576.
%
\end{thebibliography}
\end{document}